\documentclass[reqno,10pt]{amsart}
\usepackage{url}

\usepackage{graphicx}
\usepackage{multicol}
\usepackage{footmisc}

\usepackage{pgfplots}
\usepackage{subfigure}
\pgfplotsset{compat=1.12}

\usepackage{amsfonts,amsmath,amssymb}
\usepackage[latin1]{inputenc}
\usepackage{dsfont}
\usepackage{enumerate}
\usepackage{xcolor}
\usepackage[all]{xy}
\def\notshow#1\notshowend{} %
\usepackage{hyphenat}

\usepackage{graphicx}

\usepackage{tikz}
\usepackage{tikz-3dplot}
\usepackage{pgfplots}

\usepackage{multirow}
\usepackage{array}

\newcommand{\C}{\mathcal{C}}

\newcommand{\M}{\mathcal{M}} 
\newcommand{\SM}{S}

\usepackage{amssymb}
\usepackage{amsfonts}
\usepackage{mathrsfs}
\usepackage{hyperref}
\usepackage{mathtools}

\usepackage[normalem]{ulem} 

\def\bb#1\eb{\textcolor{blue}{#1}} 
\def\br#1\er{\textcolor{red}{#1}} %
\def\bm#1\em{\textcolor{purple}{#1}} %

\newcommand{\R}{\mathds R}
\newcommand{\N}{\mathds N}

\newcommand{\ho}{g_R} 
\newcommand{\vo}{v_0} 
\newcommand{\alf}{\lambda} 
\newcommand{\hh}{h} 

\newtheorem{theorem}{Theorem}
\newtheorem{prop}{Proposition}

\newtheorem{cor}{Corollary}
\newtheorem{defi}{Definition}

\newtheorem{exe}{Example}

\newtheorem{remark}{Remark}

\newcommand{\noi}{\noindent}
\newcommand{\ben}{\begin{enumerate}}
\newcommand{\een}{\end{enumerate}}
\newcommand{\bit}{\begin{itemize}}
\newcommand{\eit}{\end{itemize}}
\newcommand{\edoc}{\end{document}}

\newcommand{\A}{\mathcal{A}}

\begin{document}
%
\title[An account on  links between  Finsler and Lorentz Geometries]{An account on  links between  Finsler and Lorentz Geometries for Riemannian Geometers}

\thanks{This is a preprint of the following chapter: M. \'{A}. Javaloyes, E. Pend\'{a}s-Recondo and M. S\'{a}nchez, An Account on Links Between Finsler and Lorentz Geometries for Riemannian Geometers, New Trends in Geometric Analysis, RSME Springer Series, vol. 10, edited by A. Alarc\'{o}n, V. Palmer and C. Rosales, 2023, Springer, reproduced with permission of Springer Nature Switzerland AG. The final authenticated version is available online at: http://dx.doi.org/10.1007/978-3-031-39916-9.}

%
%
\author[M. \'{A}. Javaloyes]{Miguel \'Angel Javaloyes}\address{Departamento de Matem\'aticas, \hfill\break\indent Universidad de Murcia, \hfill\break\indent Campus de Espinardo,\hfill\break\indent 30100 Espinardo, Murcia, Spain} \email{majava@um.es}

\author[E. Pend\'{a}s-Recondo]{Enrique Pend\'{a}s-Recondo}\address{Departamento de Matem\'{a}ticas, \hfill\break\indent Universidad de Murcia, \hfill\break\indent Campus de Espinardo,\hfill\break\indent 30100 Espinardo, Murcia, Spain }
\email{e.pendasrecondo@um.es}

\author[M. S\'{a}nchez]{Miguel S\'{a}nchez}\address{Departamento de Geometr\'{\i}a y Topolog\'{\i}a, Facultad de Ciencias \& \hfill\break\indent IMAG (Centro de Excelencia Mar\'{i}a de Maeztu), \hfill\break\indent Universidad de Granada, \hfill\break\indent 18071 Granada, Spain}\email{sanchezm@ugr.es}
%
%
%

\begin{abstract}
Some links between Lorentz and Finsler geometries have been developed in the last years,  with applications even to the Riemannian case. Our purpose is to give a brief description of them, which may serve as an introduction to recent references.
As a motivating example, we start with 
Zermelo navigation problem, where its known Finslerian description permits a Lorentzian picture which  allows for a full geometric understanding of 
the original problem.  Then, we develop some issues including: (a) the accurate description of the Lorentzian causality using Finsler elements, (b) the non-singular description of some Finsler elements (such as Kropina metrics or complete extensions of Randers ones with constant flag curvature), (c) the natural relation between the Lorentzian causal boundary and the Gromov and Busemann ones in the Finsler setting, and (d) practical applications to the propagation of waves and firefronts. 

\vspace{10mm}

\end{abstract}

\maketitle              

\newpage
\tableofcontents

\section{Introduction}
Lorentz and Finsler geometries are two quite different 
extensions  of the Riemannian one, which may serve as an arena to test and,  eventually, to extend  powerful Riemannian methods.  
Typically,
Finsler and Lorentz metrics appear  when either the  anisotropies of the space or the relativistic inclusion of time lead to replace the (Riemannian) infinitesimal Euclidean scalar products by the infinitesimal models of those geometries (i.e., possibly non-reversible norms or Lorentzian scalar products, respectively). The unification of both extensions in a single Lorentz-Finsler geometry has been considered by researchers interested in certain generalizations of General Relativity, which have received a strong impulse recently \cite{U2, U1}. 
However, at a less speculative level, some links between the Finsler and Lorentz settings  appear naturally,  as well as their unification in a Lorentz-Finsler one,  enhancing the techniques and results in both geometries and multiplying their applications even at a ``real-world'' level \cite{JPS, JPS2}. 
The purpose of the present article is to give a non-technical survey about some links between both geometries. The style is adapted  to readers with a background in Riemannian Geometry and interest in Geometric Analysis, and it may motivate the study of long references such as \cite{CJS22, FHS_Memo}. 

In a nutshell, (conformal) Lorentzian Geometry can be applied at a non-relativistic level in order to describe an object or wave that propagates  at a finite maximum velocity $v_{\hbox{\tiny{max}}}$. This velocity will resemble the speed of light in Relativity, as it allows one to construct lightcones by using the velocities of objects moving at maximum speed. Then,   the (conformally invariant) relativistic stuff about lightlike directions and causality can be reinterpreted for the description. The possible   variation of $v_{\hbox{\tiny{max}}}$ with the point and time can be directly incorporated in this picture. When the maximum velocity varies not only with these elements but also with the direction,  then Finsler Geometry comes into play, and the most general description leads to  Lorentz-Finsler metrics (which are endowed with anisotropic lightcones). 

However, there are relevant cases 
where a direction-dependent function  $v_{\hbox{\tiny{max}}}$ matches  with a classical Lorentz metric. 
Noticeably, this happens for the Zermelo navigation  problem of, say, a zeppelin whose velocity is affected by the wind.  Indeed, the anisotropy of $v_{\hbox{\tiny{max}}}$ with respect to observers on Earth is due to the direction of the wind, thus, it will disappear for observers comoving with it. 

As will be stressed in \S \ref{s3}, \S \ref{s4} and \S \ref{s5}, this stimulates the  development  of notions relative to geodesics, distances and boundaries in each one of the three settings, Riemannian, Finslerian and Lorentzian, taking into account possible applications in the last two.  Moreover,  Lorentz-Finsler methods become applicable to practical purposes, such as the monitoring of the front of a wave in an anisotropic medium or a wildfire, which is emphasized in \S \ref{s6}.


Here, we begin by considering the aforementioned Zermelo navigation problem in \S \ref{s2}, which will serve as a motivation for the remainder. This  shows a first relation between three variational problems: (a) Zermelo   minimization of  the arrival time for trajectories between two points, 
(b)  minimization of a (non-symmetric) distance by
 geodesics in the class of Finsler metrics of Randers type 
  and (c) the relativistic Fermat principle for light rays in 
  the class of  standard stationary spacetimes.
 We emphasize that the spacetime viewpoint allows one to  remove two usual constraints in the Finslerian approach of Zermelo problem:  (1)  mild wind,   overcome here by using wind Finsler structures and  SSTK spacetimes in \S \ref{SSTKsec}, and (2)~time-independence, overcome by using 
 non-stationary  Lorentz metrics.

In \S \ref{s3}, geometric applications for 
relativistic spacetimes are obtained. First, we give a brief explanation about the role of Cauchy hypersurfaces for Einstein equations in \S \ref{s3a} and the meaning of causality conditions for spacetimes in \S \ref{s3b}. Then, these elements are characterized in terms of Finslerian ones, namely,  a Randers metric in the case of a standard stationary spacetime in \S \ref{s3c} and a wind Riemannian structure for SSTK ones in \S \ref{s3d}. Notice that the former spacetimes are an extension of the latter, allowing a description of settled black holes beyond their (Killing) horizons.

 In \S \ref{s4},
 we will explore the Finslerian applications of Lorentzian Geometry. In some cases, the stationary spacetime viewpoint has suggested some results related to geodesics that turn out to hold for arbitrary Finsler metrics. But it is especially interesting the case of singular Finsler metrics such as Kropina or, more generally, the wind Riemannian ones. The viewpoint of Lorentzian Geometry desingularizes the problem allowing for a better understanding of these geometries.   Moreover, the wind Riemannian structures  provide the natural full  setting  to understand  the classification of Randers metrics of constant  flag  curvature, thus, revisiting the landmark obtained in \cite{BaRoSh04}. 

In \S \ref{s5} we deepen in the geometric applications by considering boundaries. It is worth pointing out that, in the Riemannian setting, the 
Gromov and  Busemann compactification have been widely studied since the seventies, the latter in  cases such as Hadamard manifolds or CAT(0) spaces. However, their systematic development for (possibly incomplete, non-reversible) Finslerian metrics had to wait until a specific motivation came from Lorentz Geometry. In Relativity,  the causal boundary introduced by Geroch, Kronheimer and Penrose had been computed in a limited number of cases. 
The above links between  Finsler and Lorentz metrics prove that, for 
stationary spacetimes, the causal boundary corresponds to  a general type of  Busemann boundaries  of a non-reversible Finsler metric which include the (forward and backward) Cauchy ones. These boundaries, as well as their relation with the Gromov one, are now well understood and explained here.

In \S \ref{s6} we go further both in the generality of the setting and the practical applications. Focusing on the propagation of a wave whose velocities depend on the time and direction, we show that the Lorentz-Finsler setting provides a neat geometrical picture of the evolution of the wavefront. In fact, the computations reduce to solving the ODE system given by the geodesic equations of a specific Lorentz-Finsler metric. Moreover, this approach can be applied in real-world models to obtain the evolution of any physical phenomenon that satisfies Huygens' principle, wildfires being the paradigmatic example.


\section{A motivating example: Zermelo navigation problem}\label{s2}
The classical Zermelo navigation problem was proposed by the German mathematician  Ernst  Zermelo in his 1931 paper \cite{Zer31} as follows:

\begin{quote}
	In an unbounded plane where the wind distribution is given by a vector field as a function of position and time, a ship moves with constant velocity relative to the surrounding air mass. How must the ship be steered in order to come from a starting point to a given goal in the shortest time?
\end{quote}
E. Zermelo himself solved the problem using calculus of variations, reducing the problem to solve the so-called Zermelo's equation
\[\frac{d\theta}{dt}=\sin^2(\theta)\frac{\partial W}{\partial  x^1}+\sin(\theta)\cos(\theta)
\left(\frac{\partial W}{\partial x^1}-\frac{\partial W}{\partial x^2}\right)-\cos^2(\theta)\frac{\partial W}{\partial x^2},\]
where $x^1 , x^2 $ are the coordinates of $\R^2$, $\theta$ is the angle of the trajectory of the ship with the $ x^1$-axis and $W$ is the variable wind which depends on time and position (see \cite[Eq. 459.8]{Carath67}).  In the thirties of the past century this problem received the attention of some very well-known mathematicians such as Levi-Civita, Von Mises and Mani\`a \cite{Levi-C31,Mises31,Mania37} and became one  of  the classical problems in the Calculus of Variations (see \cite{Carath67}). Zermelo problem can also be solved using Optimal Control Theory (see  the classical book \cite{BrHo75} or  \cite{BBBCG19,BCGW22,BCW21,Ser09} for recent developments), but our interest will focus on more geometrical methods, namely, the use of Finsler Geometry to solve the problem. This will be possible whenever the wind is time-independent and its contribution to the velocity  does not exceed that provided by the engine.

\subsection{The case of mild  time-independent wind} 

\subsubsection{Basic Finsler setup} \label{sFinsler}
Recall that a Finsler metric in a manifold $M$ is defined as a non-negative function
$F:TM\rightarrow [0,+\infty)$, being $TM$ the tangent bundle of $M$, which is smooth away from the zero section,
positive homogeneous of degree one, namely, $F(\mu v)=\mu F(v)$ for all $v\in TM$ and $\mu>0$, and such that for every $v\in TM\setminus \bf 0$, the symmetric bilinear form $g_v:T_{\pi(v)}M\times T_{\pi(v)}M\rightarrow \R$ is positive definite, where $\pi:TM\rightarrow M$ is the canonical projection and $g_v$ is defined as
\begin{equation}\label{fundten}
	g_v(u,w)=\frac{1}{2} \frac{\partial^2}{\partial t\partial s} F(v+su+tw)^2|_{t=s=0},
\end{equation}
for all $u,w\in T_{\pi(v)}M$. A very important element of a Finsler metric is its indicatrix $\Sigma=\{v\in TM: F(v)=1\}$. Indeed, the indicatrix determines completely the Finsler metric and the positive definiteness of $g_v$ is equivalent to having an indicatrix with a positive definite 
second fundamental form in $T_{\pi(v)}M$, thus,  of positive sectional curvature, with respect to  any Euclidean scalar product in $T_{\pi(v)}M$.  So, at each tangent space $T_{\pi(v)}M$,  the indicatrix yields the hypersurface  $\Sigma \cap T_{\pi(v)}M$, which is a   positively curved  smooth sphere   enclosing  the zero vector in its  (bounded) interior  region; see for example \cite[Theorem 2.14]{JS_Pisa}.  Any hypersurface with these properties   can be regarded as the unit sphere of a (non-necessarily symmetric) norm.   
 Indeed, given  $\Sigma$, the Finsler
metric is determined as follows: for every $v\in TM\setminus \bf 0$, there exists a unique $\mu(v)>0$ such that $\mu(v) v\in\Sigma$. Then the Finsler metric is obtained as $F(v)=1/\mu(v)$. Given a Finsler manifold $(M,F)$, it is possible to define the length of any piecewise smooth curve $c :[a,b]\rightarrow M$ as 
\[\ell_F(c)=\int_a^b F(\dot c(s)) ds.\]
Observe that this length is independent of positive reparametrizations of the curve, but when one changes the orientation of $ c$, the length could change.
This length leads to a non-necessarily symmetric distance $ d_F:  M\times M\rightarrow [0,+\infty)$ defined as the infimum 
\begin{equation}\label{dF}
d_F(x_0,y_0)=\inf_{c\in C_{x_0,y_0} (M)}\ell_F(c),
\end{equation}
where $C_{x_0,y_0} (M)$ is the subset of piecewise smooth curves  between $x_0$ and $y_0$.  Then,  {\em pregeodesics} are defined as the curves that locally minimize the length functional, namely, a small enough piece of a pregeodesic has length equal to the distance between the endpoints of that piece. Moreover, a pregeodesic will be said a {\em geodesic} if in addition it is an affine reparametrization of an arc-parametrized pregeodesic;  we will assume that the domain of each geodesic is an {\em inextendible} interval $I\subset \R$, except if otherwise specified.  
A  Finsler manifold is { \em forward (resp. backward) complete } if the domain of  its  geodesics is always unbounded from above (resp. below). 
It is possible to define two types of balls for every $r>0$ and 
$x_0  \in M$:
\begin{align}
B^+_F(x_0  ,r)&=\{ y\in M: d_F(x_0,y)<r\}, \quad \text{forward ball},\label{B+}\\
B^-_F( x_0 ,r)&=\{ y\in M: d_F(y,x_0)<r\}, \quad \text{backward ball}.\label{B-}
\end{align}

\subsubsection{ Classical Finslerian solution to Zermelo problem }  Let $M$ be a manifold, $x_0, y_0\in M$ and let us try to minimize the arrival time for a moving object from $x_0$ to $y_0$ in the following situation, which is more general than Zermelo's. 
Assume that its velocity  is prescribed at every oriented direction (such a velocity can also be interpreted as the maximum permitted velocity for the object). More precisely, for each $ u\in TM\setminus\{\mathbf{0}\}$,  its oriented direction is   $[u]=\{\lambda u: \lambda>0\}$ and the prescribed velocity
$\vo([u])= \mu  u$ 
for some  $\mu>0$.  
Formally, $\vo$ becomes a section of the bundle $TM \rightarrow SM$, where $SM$ is the sphere bundle of all the oriented directions.
  Further   assume that at every $x\in M$, the set of prescribed velocities therein $\{\vo([u]): \pi(u)=x\}$  is a positively curved smooth  hypersurface  of $T_xM$,  diffeomorphic to a sphere and enclosing the zero vector. Observe that,  as explained above, this hypersurface $\Sigma$ is the indicatrix of a Finsler metric $Z$ and, by construction, $Z\circ \vo \equiv 1$.  As we will see later, this is what happens in Zermelo problem when the wind is mild.

 Assume that $[a,b]\ni s\mapsto c(s) \in M$ is any smooth curve from $x_0$ to $y_0$ with non-vanishing velocity. When the moving object follows the trajectory determined by $c$  at the prescribed velocity, then one has a time reparametrization $[0,T]\ni t \mapsto  t(s)\in [a,b]$ such that 
\begin{equation}\label{velocity}
	\vo([\dot c(s)])=\frac{d (c\circ s) }{dt}(t(s))=\dot s (t(s)) \dot c(s)=\frac{1}{\dot t(s)} \dot c(s)
\end{equation}
 and, as a consequence,  $\dot t(s) =Z(\dot c(s))$ (as $Z(\vo([\dot c(s)]))=1$). The elapsed time by the  object is then
$$T= t(b)-t(a)= \int_a^b \dot t(s) ds=\int_a^b Z(\dot c(s)) ds,$$ 
that is, it is  the length of $ c $ computed with the Finsler metric $Z$. Therefore, in order to  find the trajectories that minimize the elapsed time, we will have to find the minimizing geodesics of $Z$ from $x_0$ to $y_0$,  so that the time $T$ in the previous formula is equal to $d_Z(x_0,y_0)$ for $d_Z$  as in \eqref{dF}.

Observe that the time-independent  Zermelo problem with mild wind is a particular case of this situation.  In such case, the subset of velocities is the translation of a sphere with the wind $W$  and this translated sphere still encloses the 0 vector.  Let us denote
by   $\ho$  the Riemannian metric having as a unit sphere the velocities without wind,  and assume that the wind is never stronger than this velocity, namely, $\ho(W,W)<1$. The pair $(\ho,W)$ is usually called the  navigation data of Zermelo problem. Then the translated sphere  is determined by those $v\in TM\setminus\{\mathbf{0}\} $ such that
\begin{equation*}
	\ho \left( \frac{v}{Z(v)} -W,\frac{v}{ Z(v) }-W\right)=1.
\end{equation*}
Solving an equation of second  order,  we finally deduce that
\begin{equation}\label{zermelo}
	Z(v)=\sqrt{ \frac{1}{\alf}\ho (v,v)+\frac{1}{\alf^2}\ho (W,v)^2}-\frac{1}{\alf}\ho (W,v), \quad
	 \hbox{where} \quad \alf=1-\ho (W,W). 
\end{equation} This is a metric of Randers type, namely, 
of the form
\begin{equation}\label{randers}
F(v)=\sqrt{\tilde h(v,v)}+\tilde \omega(v),
\end{equation}
where $\tilde h$ is a Riemannian metric and $\tilde \omega$ is a one-form on $M$ such that $\|\tilde \omega\|_{\tilde h}<1$ everywhere. The condition on $\tilde \omega$ guarantees that $F$ is positive away from the zero section.  Moreover,  it turns out that it also implies that the fundamental tensor \eqref{fundten} is positive definite (see \cite[\S 11.1]{BCS} or \cite[Corollary 4.17]{JS_Pisa}). It was shown in \cite{BaRoSh04} that the family of Randers metrics is the same as Zermelo metrics in \eqref{zermelo}, namely, all the metrics as in \eqref{randers} admit an expression as in \eqref{zermelo} for some navigation data $(\ho ,W)$ which is unique.
Observe that the original Zermelo problem was formulated for a Euclidean metric $\ho$, whereas in \cite{BaRoSh04} it was generalized to Riemannian manifolds. 
A further generalization to Finsler metrics can be found in \cite[\S 3]{Sh03} (see also \cite{JV18} for Zermelo problem in pseudo-Riemannian and pseudo-Finsler metrics).
%

\subsubsection{Solution using a stationary spacetime: Fermat principle} \label{subFermatprinciple} Here, we anticipate the use of some simple notions on spacetimes and causality; the unfamiliarized reader can look  at \S \ref{s3.0} first. 
A standard stationary spacetime is a Lorentzian manifold $(\R\times M,g)$ such that
\begin{equation}\label{statmet}
	g((\tau,v),(\tau,v))=-\Lambda \tau^2+2\omega(v)\tau+g_0(v,v),
\end{equation}
where $(\tau,v)\in \R\times T_{x_0} M\equiv T_{(t_0,x_0)}(\R\times M)$, with $(t_0,x_0)\in \R\times M$,  $\Lambda:M\rightarrow \R$ smooth and positive, $\omega$ a one-form  and $g_0$ a Riemmanian metric, both on $M$.  Observe that $\partial_t$ is a Killing field of $g$ which is timelike and we will assume that $(\R\times M,g)$ is time-oriented by $\partial_t$. 

In this setting,  
consider the problem of ``traveling'' from a  point  (``event'') $(0,  x_0)$ to a vertical line $\R\ni s\rightarrow (s, y_0)$ (``stationary observer at $y_0$''), so that the increase in the $t$ coordinate is minimum. In order to model the meaning of traveling, one considers  (future-directed) timelike or lightlike curves departing from  the event but, as the latter curves will ``move faster'' (at the prescribed ``speed of light''), we will restrict ourselves to the space of (smooth, future-directed) lightlike curves from the event to the stationary observer. 
Lightlike pregeodesics present known extremization properties
which imply that light rays are the unique local minimizers of the coordinate $t$.

\begin{remark} This underlies the {\em relativistic  Fermat 
principle}, namely, the lightlike geodesics joining the event and the observer are  the critical points for the arrival time functional in the aforementioned space of lightlike curves. 
 This Fermat  principle has been present in General Relativity from the very beginning. 
The first version is due to Hermann Weyl \cite{weyl17} for static spacetimes shortly after the irruption of Einstein field equations. In the following, we will use the extended version for stationary spacetimes developed by Levi-Civita \cite{LeCi27}. One of the key points of this version is that the isometries of the Killing field allow for a reduction of the problem to a hypersurface which is not necessarily orthogonal to the stationary observers and therefore, it is not its natural restspace. As a consequence, the velocity of the light rays measured using this hypersurface is not isotropic and provides the indicatrix of a Finsler metric.
\end{remark}

As a lightlike curve $\gamma=(t,x):[a,b]\rightarrow \R\times M$ satisfies that
\begin{equation*}
	-\Lambda {\dot t}^2+2\omega(\dot x)\dot t+g_0(\dot x,\dot x)=0,
\end{equation*}
it turns out that 
\[t(b)=\int_a^b \left(\frac{1}{\Lambda}\omega(\dot x)+\sqrt{\frac{1}{\Lambda^2}\omega(\dot x)^2+\frac{1}{\Lambda} g_0(\dot x,\dot x)} \right)ds,\]
namely, the arrival time coincides with the length of the Finsler metric of Randers type
\begin{equation}\label{Fermat}
	F(v)=\frac{1}{\Lambda}\omega(v)+\sqrt{\frac{1}{\Lambda^2}\omega(v)^2+\frac{1}{\Lambda} g_0(v,v)}.
\end{equation}
Thus, using the  relativistic Fermat principle, one obtains the following characterization of lightlike geodesics  (see for example \cite[Proposition 4.1]{CJM11}).  
\begin{prop}\label{lightlikegeo}
	A lightlike curve $(t,x):[a,b]\rightarrow \R\times M$ is a geodesic of $(\R\times M,g)$ if and only if $x$
	is a geodesic of $(M,F)$ up to parametrization. Moreover, in this case, if $(t,x)$ is parametrized with the time-coordinate, then $x$ is an $F$-unit geodesic. 
	\end{prop}	
 The Fermat metric \eqref{Fermat} is very similar to the Zermelo metric \eqref{zermelo}. Indeed, both families provide the same class of metrics. 
\begin{prop}\label{equiFermatZermelo}
	The Fermat metric \eqref{Fermat} associated with a  standard stationary spacetime as in \eqref{statmet} is the Zermelo metric with navigation data $(\ho,W)$ determined by
	\begin{equation}\label{e_FZ}
\omega=-g_0(W,\cdot), \quad \ho =\frac{1}{\Lambda+\|\omega\|_0^2}g_0.
\end{equation} 
	Conversely, a Zermelo metric with navigation data $(\ho,W)$ on $M$ is the Fermat metric of a standard stationary spacetime $(\R\times M,g)$ as in \eqref{statmet} determined 
	by
	\begin{equation}\label{e_FZ2}
	g_0=\ho, \quad \omega=-g_0(W,\cdot) \quad \text{and} \quad \Lambda=1-\ho (W,W).
	\end{equation}
	Moreover, this correspondence between standard stationary spacetimes and Zermelo data is one-to-one if we assume that $\Lambda+\|\omega\|_0^2=1$, which fixes an element in the conformal class of $(\R\times M,g)$.
\end{prop}
To give an idea of the proof, observe that \eqref{e_FZ} implies 
\begin{equation}\label{lambdadevel} 1-\ho(W,W)=1-\frac{1}{\Lambda+\|\omega\|_0^2}g_0(W,W)=1-\frac{\|\omega\|_0^2}{\Lambda+\|\omega\|_0^2}=\frac{\Lambda}{\Lambda+\|\omega\|_0^2},
\end{equation}
and putting all this together, one gets straightforwardly that the Zermelo metric with data given in \eqref{e_FZ} is just the Fermat metric in \eqref{Fermat} 
(see \cite[Proposition 3.1]{BiJa11} and \cite{GHWW} for the equivalence between Fermat and Zermelo metrics). 
Observe also that the Fermat metric is conformally invariant, namely, if we multiply $g$ by a positive function $\phi$ on $M$, then the associated Fermat metric is preserved. Using this fact, it is possible to choose an element of the conformal class of $g$ normalized in some sense. Indeed, given $(\R\times M,\tilde g=\phi g)$ with $\phi=1/(\Lambda+\|\omega\|_0^2)$, its Fermat metric satisfies that $\tilde \Lambda+\|\tilde\omega\|_0^2=1$, where $\tilde\Lambda$, $\tilde \omega$ and $\tilde g_0$ are the data of $\tilde g$ and the norm is computed with $\tilde g_0$. In this case, it follows from \eqref{lambdadevel} that $\tilde\Lambda \coloneqq \phi\Lambda=1-\ho(W,W)$. 
\begin{remark}\label{lightZermelo}
	Proposition \ref{equiFermatZermelo} makes apparent the following non-relativistic interpretation of the Fermat metric. The indicatrix of the Fermat metric is the subset of the velocities of light rays measured by the stationary observers in the tangent space to the slices $\{t_0\}\times M$, which is the translation of a Riemannian sphere.  Clearly, this interpretation is not relativistic (all relativistic observers measure the same {\em speed of light}!). Indeed, a relativistic observer would use their space at rest, which is orthogonal to their timelike direction.  This space is infinitesimal and can be identified with the tangent space to the slices only when $\omega=0$. 

	\end{remark}
 It is worth pointing out that the stationary spacetime associated with a Zermelo metric can be interpreted as a sort of Analogue Gravity (see  \cite[\S 2]{BLV05} and \cite{DS19}). 

\subsection{The case of arbitrary time-independent wind}
\subsubsection{Emergence of wind Riemannian structures}
\label{subsubsec:wind_Riem}
Let us assume now that in the navigation data $(\ho,W)$ the wind $W$ is arbitrary, namely, let us remove the constraint $\ho(W,W)<1$. This means that we will translate the spheres of $\ho$ with a vector which is possibly not contained in the sphere.  The set of all the  translated spheres will form a smooth hypersurface $\Sigma\subset TM$ which generalizes the  indicatrix  
of a Finsler metric and is called a {\em wind Riemannian structure}. 
 
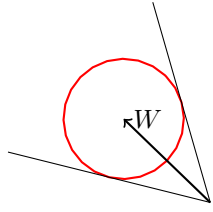
\begin{figure}
\centering
	\begin{tikzpicture}[scale=0.8,rotate around={46:(0,-2)}]
	\draw [red,thick,domain=0:360] plot ({cos(\x)}, {sin(\x)});
	\draw (0,-2) -- (1.732,1);
	\draw (0,-2) -- (-1.732,1);
	\draw[->,thick] (0,-2) -- (0,0) node[below,right]{$W$};
	\end{tikzpicture}\hspace{1cm}
	\caption{In the presence of a strong wind $W$ with $\ho(W,W)>1$, the only directions allowed for the moving object are those which intersect the red sphere.}
\label{fig:strong_wind}
\end{figure}

 So, each translated sphere $\Sigma_x$  plays the role of an indicatrix at $x\in M$, which might  not enclose the zero in its interior region; thus the regions with $\ho(W,W)>1$, $\ho(W,W)=1$ or $\ho(W,W)<1$ are called of  {\em strong, critical}  or  {\em mild wind}, respectively.  In particular, when the wind is   strong,  
 one has some {\em admissible directions} at each $x\in M$, defined as the oriented directions $[u]$ which intersect the indicatrix; they form a closed solid cone in $T_xM$. Those $[u]$ in the interior of the cone intersect twice the indicatrix and, so, provide naturally two velocities $v_+([u]), v_-([u])$ (the latter contained in the segment between the origin and the former); both velocities merge in a single one for the admissible directions in the boundary (see Figure \ref{fig:strong_wind}). 

 In the region of strong wind, Zermelo problem for $x_0, y_0\in M$ splits into several ones. First,  determine whether there exists an admissible curve from $x_0$ to $y_0$. 
  If this is the case, prescribe either $v_+$ or $v_-$ for the moving object and find extremal curves for the arrival time. In the  case of $v_+$, it is natural to wonder for a (local or global) minimizer of the arrival time and one properly has the previous Zermelo problem. Indeed, such a problem  makes sense for curves which may cross regions of strong, critical  and mild wind, and $v_+$ will match smoothly with $\vo$ on the whole $M$.  In the case of prescribing $v_-$, the natural problem would be to find  {\em maximizers} of the arrival time entirely contained in the region of strong wind. This  case would occur when the object  tries to delay the arrival at $y_0$ as much as possible by making its engine power go against the direction of the wind. Technically, this makes sense because the portion of the indicatrix corresponding to $v_-$ is concave (as happens for the unit timelike directions of a Lorentz metric). 


 These problems were first studied in \cite{CJS22}. Formally, previous ideas of the mild case work similarly 
to obtain a Finsler metric which measures the elapsed time along an admissible trajectory.  In the strong wind case we obtain two metrics:
\begin{align}\label{zermeloZ}
	Z(v)=&\frac{1}{\alf}(\sqrt{ \alf \ho(v,v)+\ho(W,v)^2}-\ho(W,v)),\\
	\label{zermeloZl}  Z_l(v)=&\frac{1}{\alf}(-\sqrt{ \alf \ho(v,v)+\ho(W,v)^2}-\ho(W,v)),
\end{align} 
being defined both  of them  in 
\[   A_l =\{ v\in TM: \alf<0, \ \alf \ho(v,v)+\ho(W,v)^2>0, \ -\ho(W,v)<0\}.\] 
Observe that if we want to compute the minimizing solutions of Zermelo problem, then we have to consider $Z$,  which corresponds to $v_+$ above.  
The    above expressions do not hold in the critical case  $\ho(W,W)=1$, but it is possible to obtain expressions 
valid for arbitrary winds  multiplying  by the conjugate  in \eqref{zermeloZ}, \eqref{zermeloZl}, i.e.,  
\begin{align}\label{zermeloArbitrary}
	Z(v)=&\frac{\ho(v,v)}{\sqrt{\alf \ho(v,v)+\ho(W,v)^2}+\ho(W,v)},\\
	\label{ZermeloArbitraryZl} Z_l(v)=&\frac{-\ho(v,v)}{\sqrt{\alf \ho(v,v)+\ho(W,v)^2}-\ho(W,v)}.
\end{align} 
Observe that now $Z$ is defined for an arbitrary $W$ in the domain 
\[A=\{ v\in TM \setminus \mathbf{0}: \text{$\lambda>0$ or }  \alf \ho(v,v)+\ho(W,v)^2>0, \ \ho(W,v)>0\}\]
and, in the case of critical wind, the metric $Z$ is of Kropina type, namely, the quotient of a Riemannian metric by a one-form:\footnote{These metrics are well-known in the Finsler literature since the original article \cite{Kr}.}
\[Z(v)=\frac{\ho(v,v)}{2 \ho(W,v)}.\]
In particular, the domain $A$ coincides with $TM\setminus \bf 0$ when $\alf>0$  (the region of mild wind), it is the half-space $\ho(W,v)>0$ in the tangent bundle of the region 
\[M_{crit}=\{x\in M: \alf=0\}\]
(the region of critical wind) and a conic region in 
\[M_l=\{x\in M: \alf<0\}\] 
(the region of strong wind). On $M_l$, the domain $A$ coincides with the timelike vectors of the Lorentzian metric $-\hh$ on $M_l$,\footnote{Notice that, whenever $ \lambda < 0 $, $\hh$ has signature  $(+,-,\dots , -)$. Anyway, $\hh$ is well defined as a signature changing metric on the whole $M$, so that it can be used to determine the admissible curves between any two points $x_0, y_0\in M$. Comparing with $\tilde h$ in \eqref{randers} (which was defined for $\alf>0$), one has $\tilde h=h/\alf^2$.} where 

\begin{equation}\label{-tildeh}
 \hh(v,v)=   \alf \ho(v,v)+\ho(W,v)^2
\end{equation}
in the half-space $h(W,v)  <0$. 
Moreover, the domain of $Z_l$ is $A_l=A\cap TM_l$. 
It turns out that in the region $M_{crit}\cup M_l$, the metric $Z$ is conic, namely, at every $x\in M_{crit}\cup M_l$ it is not defined in the whole $TM$ but in  the conic region $A_x=A\cap T_xM$. Moreover, it is positive homogeneous of degree one and its fundamental tensor \eqref{fundten} is positive definite. The metric $Z$ can be extended to the boundary of $A_x$, but this extension is not smooth and its fundamental tensor cannot be extended to the boundary. However, $Z$ is a classical Finsler metric in the region $M\setminus \{M_{crit}\cup M_l\}$. On the other hand, the metric $Z_l$ is always  {\em conic}\footnote{See \cite{JS_Pisa} for a full development of this condition.} and its fundamental tensor has index $n-1$. Moreover, it can be extended to the boundary of $A_x$, this extension coincides with that of $Z$, but again its fundamental tensor does not admit such extension (see \cite[\S 3.3 to \S 3.5]{CJS22}).

 Summing up, Zermelo problem is modeled  in terms of a wind Riemannian structure, whose  Finslerian description retains some elements of the mild wind case. However, it also includes new ingredients (Kropina metric, $F_l$ with  concave indicatrix) which become complicated and apparently singular. Next, the spacetime viewpoint will simplify the picture giving an elegant solution. 

\subsubsection{Solution using an SSTK spacetime}\label{SSTKsec}
The process to solve Zermelo problem for mild wind using a spacetime developed in \S \ref{subFermatprinciple} can be extended for a general wind just by considering Lorentzian metrics as in (9) with an arbitrary $\Lambda$ (not necessarily positive)  satisfying  $\Lambda +\|\omega\|_0^2>0$;  this condition must be imposed because it  is equivalent to the Lorentzian character of \eqref{statmet}.  Following \cite{CJS22},  this class of spacetimes is called  {\em SSTK (standard with a space-transverse Killing vector field)} spacetimes. Observe that the points $x\in M$ with $\Lambda(x)<0$ correspond to vertical lines $s\rightarrow (s,x)\in\R\times M$ which are spacelike in $(\R\times M,g)$, and then $\partial_t$ is not timelike  therein. It is still consistent though to consider  the time-orientation determined by $dt>0$, namely, a timelike vector $(\tau,v)$ is future-directed if and only if $\tau>0$.  It is worth pointing out 
that the usual relativistic Fermat principle does not apply to such a vertical line,  as it is not necessarily timelike; however, as shown in \cite[Theorem 7.4]{CJS22}, this principle can be extended for arrival curves of arbitrary causal character.  Thus, the solution of Zermelo problem in this setting can be described as follows: 
\begin{quote}
Given $x_0,y_0\in M, x_0\neq y_0$ and the Zermelo data $(\ho,W)$, construct an SSTK spacetime with $g_0, \omega$ and $\Lambda$ given as in  \eqref{e_FZ2}.  An admissible curve $c:[0,T]\rightarrow M$ from $x_0$ to $y_0$ with prescribed velocity in the wind Riemannian structure $\Sigma$ is a critical trajectory for the arrival time $T$ if and only if the future-directed lightlike curve $[0,T]\ni t\mapsto (t,c(t))\in \R\times M$ 
is a lightlike pregeodesic of the spacetime. 

 In this case, the arrival time $T$ of $c$ is globally minimizing (resp. maximizing) if and only if   $(T,y_0)\in \R\times M$   is the first (resp. last) $t$ of the observer $\{(t,y_0), t\in \R\}$ at $y_0$
reached by future-directed causal curves starting at $(0,x_0)$. 
\end{quote}
 We can also recover $c$ as a unit speed geodesic  for $F$ or $F_l$ (depending whether it is locally minimizing or maximizing) by repeating the process of \S \ref{subFermatprinciple}. Indeed,  it follows that a lightlike curve 
$\gamma=(t,x):[a,b]\rightarrow \R\times M$ is a lightlike pregeodesic if and only  if its projection $x$ is either a pregeodesic of one of the conic metrics 
\begin{align}\label{fermatArbitraryF}
	F(v)=&\frac{g_0(v,v)}{-\omega(v)+\sqrt{\omega(v)^2+\Lambda g_0(v,v)}},\\
	\label{fermatArbitraryFl} F_l(v)=&\frac{-g_0(v,v)}{\omega(v)+\sqrt{\omega(v)^2+\Lambda g_0(v,v)}},
\end{align} 
 or it is a suitable curve with velocity in the  closure of the conic domain  constructed from the lightlike pregeodesics of $\hh$  \cite[Theorem 5.5]{CJS22}. 

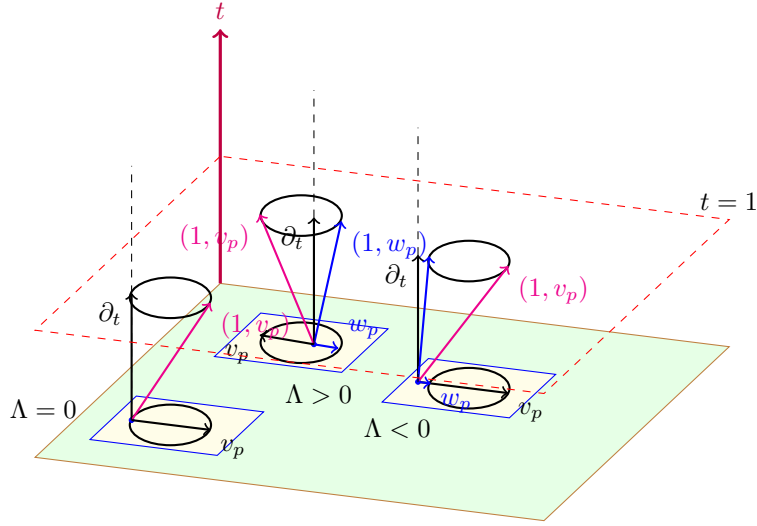
\begin{figure}
\centering
\tdplotsetmaincoords{70}{110}
\begin{tikzpicture}[scale=1.8,tdplot_main_coords]
    \filldraw[
        draw=brown,%
        fill=green!10,%
    ]          (0,0,0)
            -- (4,0,0)
            -- (4,4,0)
            -- (0,4,0)
            -- cycle;
       \filldraw[
        draw=blue,%
        fill=yellow!10,%
    ]          (1,2,0)
            -- (2,2,0)
            -- (2,3,0)
            -- (1,3,0)
            -- cycle;
          \filldraw[
        draw=blue,%
        fill=yellow!10,%
    ]          (0.5,0.5,0)
            -- (1.5,0.5,0)
            -- (1.5,1.5,0)
            -- (0.5,1.5,0)
            -- cycle;
            \filldraw[
        draw=blue,%
        fill=yellow!10,%
    ]          (2.5,0.25,0)
            -- (2.5,1.25,0)
            -- (3.5,1.25,0)
            -- (3.5,0.25,0)
            -- cycle;
           \draw[thick] (1,1,0) circle (0.3cm and 0.15cm);
           \filldraw[blue] (1,1.1,0) circle (0.5pt);
           \draw[thick] (1,1,1) circle (0.3cm and 0.15cm);
           \draw[thick,->] (1,1.1,0) -- (1,1.1,1) node[left,anchor=north east]{$\partial_t$};
           \draw[thick,blue,->] (1,1.1,0)  -- (1+0.01,1.1+0.22,1) node[left,anchor=north west]{$(1,w_p)$};;
           \draw[thick,magenta,->] (1,1.1,0)  -- (1-0.05,1.1-0.445,1-0.05) node[left,anchor=north east]{$(1,v_p)$};
           \draw[thick,->] (1,1.1,0)  -- (1-0.05,1.1-0.445,0) node[left,anchor=north east]{$v_p$};
            \draw[thick,->,blue] (1,1.1,0)  -- (1+0.02,1.1+0.2,0) node[right,anchor=south west]{$w_p$};
           	\node at (2.5,2.3,0)(s){$\Lambda<0$};
           \draw[thick] (3,0.7,0) circle (0.3cm and 0.15cm);
           \filldraw[blue] (3,0.395,0) circle (0.5pt);
           \draw[thick] (3,0.7,1) circle (0.3cm and 0.15cm);
           \draw[thick,->] (3,0.395,0) -- (3,0.395,1) node[left,anchor=north east]{$\partial_t$};
           \draw[thick,magenta,->] (3,0.395,0)  -- (3,0.395+0.62,1) node[left,anchor=north west]{$(1,v_p)$};
           \draw[thick,->] (3,0.395,0)  -- (3,0.395+0.62,0) node[left,anchor=north west]{$v_p$};
           	\node at (3,-0.3,0)(s){$\Lambda=0$};
           \draw[thick] (1.5,2.5,0) circle (0.3cm and 0.15cm);
           \filldraw[blue] (1.5,2.1,0) circle (0.5pt);
           \draw[thick,->] (1.5,2.1,0) --(1.5,2.1,1) node[left,anchor=north east]{$\partial_t$};
           \draw[thick] (1.5,2.5,1) circle (0.3cm and 0.15cm);
           \draw[thick,magenta,->] (1.5,2.1,0) -- ((1.5,2.1+0.71,1) node[left,anchor=north west]{$(1,v_p)$};
           \draw[thick,->] (1.5,2.1,0) -- ((1.5,2.1+0.71,0) node[left,anchor=north west]{$v_p$};
           \draw[thick,->,blue] (1.5,2.1,0) -- ((1.5,2.1+0.1,0) node[left,anchor=north west]{$w_p$};
           	\node at (2,1.5,0)(s){$\Lambda>0$};
           \draw[thick,blue,->] (1.5,2.1,0)  -- (1.5,2.1+0.088,1);
            \draw[dashed] (1,1.1,0) -- (1,1.1,2);
             \draw[dashed] (3,0.395,0)  -- (3,0.395,2) ;
              \draw[dashed] (1.5,2.1,0) -- (1.5,2.1,2);
    \draw[very thick,->,purple] (0,0,0) -- (0,0,2) node[left,anchor=south]{$t$};
     \draw[
        draw=red,dashed%
    ]          (0,0,1)
            -- (4,0,1)
            -- (4,4,1)
            -- (0,4,1) node[left,anchor=south]{$t=1$} 
            -- cycle;
\end{tikzpicture}
\caption{We show how to obtain the wind Riemannian structure associated with the SSTK spacetime by intersecting the lightlicone with the slice $t=1$ in the three different cases: $\Lambda=0$ (critical wind), $\Lambda>0$ (mild wind) and $\Lambda<0$ (strong wind). Observe that in the first two cases $v_p$ is the projection of a unique lightlike vector, while in the last one there are exactly two lightlike vectors projecting and pointing out in the same direction as $v_p$.}
\label{figurecones}
\end{figure}

As in the previous subsection,  now the Fermat metric \eqref{Fermat} splits into two, and one has to multiply by the conjugate to obtain the  metrics \eqref{fermatArbitraryF} (defined on the whole $M$) and \eqref{fermatArbitraryFl} (only on $M_l$); compare with those 
in \eqref{zermeloZ}, \eqref{zermeloZl}, \eqref{zermeloArbitrary} and \eqref{ZermeloArbitraryZl}. This splitting is naturally interpreted now,  because there are two future-directed lightlike vectors which project  onto each $v\in TM_l$ (see Figure \ref{figurecones}).  Extending Proposition \ref{equiFermatZermelo}, 
 these two metrics can be identified respectively with those  in \eqref{zermeloArbitrary} and \eqref{ZermeloArbitraryZl}.   
  

 Summing up, the spacetime viewpoint yields a full solution of Zermelo problem which permits a unified description of the solutions and allows one to recover the (apparently singular) Finslerian description in the spacelike part. In particular, the   lightlike geodesics in the spacetime, which solve the problem, can be characterized as follows \cite[Theorem 5.5]{CJS22}:

\begin{theorem}\label{lightlikegeo2}
	A lightlike curve $(t,x):[a,b]\rightarrow \R\times M$ is a geodesic of $(\R\times M,g)$ if and only if 
it lies in one of the following three cases: (a) $\dot x$ is entirely contained in 
\begin{equation}\label{conicA}
A=\{ v\in TM\setminus \mathbf{0}: \; \hbox{either} \; \Lambda>0 \; \hbox{or} \;  \Lambda g_0(v,v)+\omega(v)^2>0, \; \omega(v)<0\}
\end{equation}
 and, in this case, $x$ is a pregeodesic of either $F$ or $F_l$, (b) $x$ is constantly equal to some $x_0$ with $d\Lambda(x_0)=0$, or (c) $x$ is a suitable curve with $\dot x$ in the closure of $A$.\footnote{It is worth emphasizing that this last case appears in the same footing as the others from the spacetime viewpoint. However, in the Finslerian description, it corresponds with the limit of the geodesics of $F$ and $F_l$ (indeed, as a limit, $F(\dot x)=F_l(\dot x)\equiv 1$). Moreover, $x$ might start at the region of $M_l$, arrive at the region of critical wind $M_{crit}$ and come back to $M_l$. On $M_l$,  $x$ becomes a pregeodesic of $\Lambda g_0+ \omega\otimes \omega$ consistent with \eqref{-tildeh} (see the case (iii) (b) in the aforementioned \cite[Theorem 5.5]{CJS22}).}
	

	\end{theorem} 
\subsection{The time-dependent case}\label{sec:time-dep}
 The  Zermelo problem in a manifold $M$ with a time-dependent wind  can be handled by  using time-dependent Finsler metrics.  
This was  done  by Mani\`a \cite{Mania37} as well as by Markvorsen using frozen metrics \cite{M17}. Here, we will develop  a spacetime picture as a natural extension of our  framework. Indeed,  we
 will consider again a more general problem, namely,  we will assume that the velocity is prescribed at every direction 
 but,  now, this prescription may  have a dependence on time.  Thus, one has a smooth hypersurface $\Sigma_{(t_0,x_0)}$ of $T_{x_0}M$ diffeomorphic to a sphere and positively  curved   
 at each $t_0\in \R$, and all these hypersurfaces vary smoothly with $(t_0,x_0)$, providing a smooth submanifold\footnote{See \cite[Definition 2.8 and Example 2.16]{CJS22} for a subtlety about  smoothness  applicable here.} of $T(\R\times M)$  with codimension 2.

\begin{figure}
	\centering
	
	\begin{tikzpicture}[scale=0.6,rotate around={-46:(-4,0)}]
	\draw[blue] 
	(-6.1,4.5) -- (-4,0) -- (-1.9,4.5);
	\draw[thick,color=purple] (-5.3,2.45) -- (-6.2,2) -- (-6.2,0.2) -- (-1.7,2.3) -- (-1.7,4) -- (-2.2,3.77);
	\draw[blue] (-4,4) ellipse (1.83cm and 0.5cm);
	\draw[purple,rotate around={26:(-3.88,2)}] (-3.88,2) ellipse (0.943cm and 0.5cm);
	\node at (-4,-0.3)(s){$p$};
	\node at (-1.4,4.5)(s){$\C_p$};
	\node at (-1.2,3.5)(s){$\Pi$};
	\node at (-3.4,2.95)(s){$\Sigma$};
	
	
	\end{tikzpicture}
	\caption{The cone structure $\C_p$ is intersected with an affine hyperplane $\Pi$ in a strongly convex hypersurface $\Sigma$ of $\Pi$.}
\label{cone} 
\end{figure}
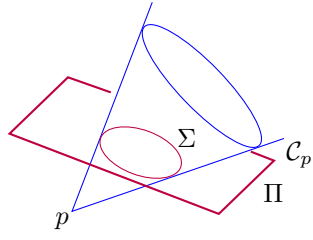

  Observe that in the non-relativistic spacetime $\R\times M$, being the coordinate $\R$ the (absolute)  time, a curve
   $s\rightarrow (t(s),x(s))$ corresponds to a curve with velocity $v(s)=\frac{\dot x(s)}{\dot t(s)}\in T_{x(s)}M$ at the instant $s$ (recall \eqref{velocity}). This means that if the velocities 
    at each time $t_0$  are prescribed by the hypersurface $\Sigma_{t_0}  =\cup_{x_0\in M} \Sigma_{(t_0,x_0)}\subset 
   TM$,  
   then the  vectors
 $(\tau,v)\in T_{(t_0,x_0)}(\R\times M)$  
     tangent  to the allowed curves (traveling at the prescribed velocities) must satisfy that $v/\tau\in \Sigma_{t_0}$ (with $\tau>0$). It turns out that  these  allowed tangent vectors 
determine a cone structure $\C$ (see \cite[Definition 2.7]{JS20}), namely, the smooth embedded hypersurface $\C \subset T(\R\times M)\setminus 
\mathbf{0}$ given by
$$
\C_{(t_0,x_0)} = \{\tau (1,v)\in T_{(t_0,x_0)}(\R\times M): \tau>0, \ v\in \Sigma_{(t_0,x_0)}\},
$$
where $\C_{(t_0,x_0)}=\C\cap T_{(t_0,x_0)}(\R\times M)$ for every $(t_0,x_0)\in \R\times M$. The requirements on each $\Sigma_{(t_0,x_0)}$ become equivalent to the {\em strong convexity} of $\C$ \cite[Proposition 2.26]{JS20}, that is, when $\C_{(t_0,x_0)}$    is intersected by an affine hyperplane of $H\subset T_{(t_0,x_0)}(\R\times M)$ which is not  tangent to $\C_{(t_0,x_0)}$, then $H\cap \C_{(t_0,x_0)}$  is a strongly convex hypersurface in the hyperplane (in fact, a positively curved sphere); see Figure \ref{cone}.

 When the natural vector field $\partial_t$ is timelike for $\C$ (i.e., it lies in the interior of the solid cone bounded by $\C$ at each point) the 
 easiest way to describe  $\C$ is using a Lorentz-Finsler metric $G$ in $\R\times M$ constructed as follows.  Let $F_{t_0}$ be the Finsler metric in 
 $M$ whose indicatrix coincides with $\Sigma_{t_0}$, that is,
$$\Sigma_{t_0}=\{v\in TM: F_{t_0}(v)=1\}.$$ 
  Then, define $G$   at every point $(t_0,x_0)$ as
\begin{equation}\label{metricG}
G(\tau,v)=\tau^2-F_{t_0}(v)^2, \quad  \forall (\tau,v)\in T_{(t_0,x_0)}(\R\times M). 
\end{equation}
\ Observe that $G$ is a  Lorentz-Finsler metric\footnote{In the sense of \cite[\S 3.2]{JS20} (thus, positive and two-homogeneous) up to the fact that it may be non-smooth on $\partial_t$. This is not relevant here, because we are only interested in the directions of $\C$ (anyway, $G$ can be smoothen on $\partial_t$; see \cite[Theorem 5.6]{JS20}.)} with fundamental tensor of  Lorentzian type (replace $F^2$ with $G$ in \eqref{fundten}). It   turns out that the cone structure $\C$ is given by the lightlike vectors of $G$  (with $\tau>0$).  Indeed,    $G(\tau,v)=0$ if and only if $\tau^2=F_{t_0}(v)^2$, which is equivalent to 
$v/\tau\in \Sigma_{t_0}$,  whenever  $\tau>0$. 

 Once $G$ is obtained,  the admissible trajectories between an event $(0,x_0)$ and an  ``observer''  $\R\ni t\rightarrow (t,y_0)$ at a different  point $y_0\in M$ are  lightlike curves of $(\R\times M,G)$.   The Fermat principle can also be extended to the Finslerian relativistic case (see, \cite{Perlick06}). Then,  it follows that the solutions of the time-dependent case are provided by the projections  $x$ on $M$ of the lightlike pregeodesics $(t,x)$ of $(\R\times M,G)$ 
going from the event to the observer, 
with $\dot t>0$. 

 Finally,  the case when $\partial_t$ is not timelike can be handled by using an auxiliary timelike vector field $X$ which would represent comoving observers; see \cite[Remark 6.1]{JPS}. 


\section{Relativistic applications} \label{s3}

\subsection{Basic Lorentz setup}\label{s3.0} Recall first some basic notions related to relativistic spacetimes following  \cite{BEE,O}, 
which are standard for readers with background in Riemannian Geometry.
Let $g$ be a  Lorentz metric on an $ (n+1) $ manifold  $\M$  with  signature $(-,+, \dots , +$), which implies that each connected component of $\M$ will be either non-compact or of zero Euler characteristic. A non-zero tangent vector $v\in T\M$ is called {\em timelike}, {\em lightlike} or {\em spacelike} if  $g(v,v)<0, g(v,v)=0$ or $g(v,v)>0$, respectively; {\em  causal} vectors are the timelike or lightlike ones and {\em null} vectors the lightlike or zero ones.\footnote{The convention for the zero vector depends on the reference, we follow \cite{MS}.} It is obvious that any conformal Lorentz metric $g^*=\Omega g$ (for some function $\Omega$) will have the same cones as $g$ and, remarkably, the converse is true, that is, the lightlike vectors determine the conformal class of the metric (see for example \cite[Lemma 2.1]{BEE}).

$(\M,g)$ is called a {\em spacetime} when $\M$ is connected and one assumes that a time-orientation has been chosen. The latter is a continuous choice of one of the two  cones determined by the causal vectors at each tangent space $T_p\M, p\in \M$; each chosen cone will be called {\em future}, while  the non-chosen one will be {\em past}. The properties of time-orientations are similar  to those of usual orientations on $\M$ (including the existence of a time-orientable double covering), but both notions are independent. In spacetimes, the names of {\em future-} or {\em past-directed} are applied directly to causal vectors; moreover, all the related names (timelike, future or past-directed etc) are transferred directly to (smooth) curves
taking into account the character of the velocity. Submanifolds are called spacelike, timelike or lightlike depending on whether the induced metric has index 0, 1 or it is degenerate, respectively.

A major topic of research in Mathematical Relativity with interest in Geometric Analysis  is the {\em initial value problem} for Einstein equations, which is posed on {\em Cauchy hypersurfaces}, 
see  for example \cite{C-B_libro,ChF,N,R}. Next, we will describe briefly the setting and, then, some applications of the aforementioned geometric links with Finsler manifolds.

\subsection{The initial value problem and Cauchy hypersurfaces}\label{s3a}

The Einstein equation for a  spacetime $(\M ,g)$,  $n+1>2$ , is
\begin{equation}\label{e_Einstein}
\hbox{Ric} + \frac{1}{2}\, \hbox{Scal} \, g +\Lambda_c g= T,
\end{equation}
where Ric and Scal denote, respectively, the Ricci tensor and scalar curvature of $g$,   $\Lambda_c\in \R$ is the {\em cosmological constant} and $T$ is the {\em stress-energy tensor} up to a constant (which depends on units and provides the suitable coupling between matter and geometry). Roughly speaking, the initial value problem will consist in  specifying initial data $(h, \A)$ 
on some 3-manifold $\SM$ and finding a spacetime $(\M,g)$ satisfying \eqref{e_Einstein}, where $\SM$ is  embedded as a spacelike hypersurface with metric $h$ and second fundamental form $\A$. In particular, the initial data must satisfy the constraints associated with the Gauss and Codazzi equations, which turn out essential for the analysis of solutions. In general, one should also specify initial data on $\SM$ for $T$ as well as additional equations to be satisfied for the specific type of matter modeled by $T$ (so that \eqref{e_Einstein} remains as a hyperbolic system of equations in $g$ and $T$). However,  the latter are not required for the vacuum case $T=0$ (which may be a good approximation for modeling the empty space outside a star).
It is worth pointing out that,  when $T=0$, (semi-)Riemannian  Schur's lemma reduces \eqref{e_Einstein} to the equation of Einstein manifolds Ric$= a g, a\in \R$ and, in the case of $\Lambda_c=0$, to the Ricci flat case Ric $=0$.

Global existence and uniqueness of solutions for the vacuum case and other $T$'s have been proved,  starting at the classical results by Choquet-Bruhat and Geroch \cite{C-B,C-B_G,G}.
A posteriori, the  initial data manifold $\SM\subset \M$ will be a (smooth, spacelike) {\em Cauchy hypersurface} of the solution $(\M,g)$, i.e., $\SM$ is crossed exactly once by any inextendible causal curve.\footnote{Sometimes Cauchy hypersurfaces and the initial value problem are allowed to be more general so that non-smooth data or data on degenerate hypersurfaces are permitted, but we will not go into these issues.} The uniqueness of the crossing point is related to the existence of solutions (otherwise, data on a point of $S$ might influence on data on points of $S$ reached by a future-directed causal curve), and  the 
existence  of the crossing point is related to the uniqueness of the solutions.  Indeed,  the existence of a {\em unique} maximal solution is ensured {\em assuming that $\SM$ remains Cauchy}, but different extensions might exist if this assumption is dropped (such a possibility is related to the so-called {\em strong cosmic censorship hypothesis}; see for example \cite[\S 12.1]{Wald}).
 
Summing up, a posteriori, the solution $(\M,g)$ to the initial value problem of the Einstein equation will admit a Cauchy hypersurface $\SM$ and $(\M,g)$ can be regarded as the domain of dependence of the data on $\SM$. Thus,  it arises the natural question of determining 
 the class of spacetimes admitting such a hypersurface as well as to determine when a given spacelike hypersurface is Cauchy, which will be adressed in the next two subsections.

\subsection{Global hyperbolicity and causality of  spacetimes}\label{s3b}

Recall first the following basic elements of causality for a spacetime $(\M, g)$ (see \cite{MS} for a complete study). 
For $p,q\in \M$, we say that $p$ lies in the chronological (resp. strict causal) past of $q$ if there exists a future-directed timelike (resp. causal) curve from $p$ to $q$; in this case, we write $p\ll q$ (resp. $p<q$). The chronological (resp. causal) future of $p$ is defined as
\begin{equation}
\label{e_causalrelations} 
I^+(p)=\{q\in\M :  p\ll q\} \quad (\hbox{resp.}\; J^+(p)=\{q\in\M :  p< q\}\cup \{p\}),
\end{equation}
then, one defines analogously the chronological and causal pasts $I^-(p), J^-(p)$. 
Following \cite{BS07}, a spacetime is called {\em globally hyperbolic} when 
\begin{equation}\label{e_def_glob_hyp}
p\not< p \quad \hbox{and}\quad J^+(p)\cap J^-(q) \; \hbox{is compact} 
\end{equation}
for all $p,q\in\M$. The first condition means that the spacetime does not contain  any causal loop, that is, $(\M,g)$ is {\em causal}, and the second one that no observer can see the sudden appearance or disappearance of a particle from the spacetime, i.e., $(\M,g)$ does not contain {\em naked singularities}.
Taking into account Geroch's \cite{G} and the further developments in \cite{BS05,BS06}, the following  characterization holds.
\begin{theorem}\label{t_Cauchytemporal}
A spacetime is globally hyperbolic if and only if it admits a (smooth, spacelike) Cauchy hypersurface. In this case, there exists some onto smooth map  $\tau: \M \rightarrow \R$  such that the spacetime splits as an orthogonal product
$$
\M \cong \R\times M, \quad g=-\Lambda d\tau^2+g_t, 
$$ 
where $\Lambda>0$ is a function on $\M$, $g_t$ is a Riemannian metric on each slice $\{t\}\times M$  varying smoothly with $t$, and  each slice 
becomes a Cauchy hypersurface in $\M$.\footnote{The natural time-orientation chosen here and in the remainder is the one such that the timelike vector $\partial_\tau$ becomes future-pointing.} Such a $\tau$ is called a {\em Cauchy temporal function}. Moreover,  any Cauchy hypersurface $S$ can be regarded as a level of some Cauchy temporal function, obtaining then an orthogonal splitting as above such that $S \cong \{0\}\times M$.
\end{theorem}
A generalization of globally hyperbolic spacetimes are the {\em causally simple} ones, where the second hypothesis in \eqref{e_def_glob_hyp} is weakened into the condition: $J^+(p), J^-(p)$ are closed for all $p,q$. For these spacetimes, no splitting as in Theorem \ref{t_Cauchytemporal} can be ensured, however, they  are {\em stably causal}, that is, they admit a {\em temporal function} $\tau$. This means that $\tau: \M\rightarrow \R$ is onto  with grad$(\tau )$ timelike and past-directed (thus, $\tau$ grows strictly on any future-directed causal curve and the  slices $\tau=$ constant are spacelike).  
A simple example of these properties is the following
(see \cite[Theorem 3.67]{BEE} and Theorem \ref{t_escala_causal} below).

\begin{prop}\label{p_escala_causal} Let $(M,g_R)$ be a Riemannian manifold and consider the natural product spacetime $(\M=\R\times M,g=-d\tau^2+g_R)$. Then, the canonical projection $\tau: \R\times M\rightarrow \R$ is a temporal function and

\smallskip


\noi (A) $(\M, g)$ is causally simple if and only if $g_R$ is convex (i.e., each $p,q\in M$ can be joined by a minimizing geodesic).

\smallskip
\noi
 (B) The following properties 
are equivalent:
\begin{enumerate}
\item[(B1)] $(\M, g)$ is globally hyperbolic.

\item[(B2)] The slices $\{t\}\times M$ are Cauchy hypersurfaces (i.e., $\tau$ is Cauchy temporal).

\item[(B3)] $g_R$ is complete.
\end{enumerate}
\end{prop}

\begin{remark}\label{r_static} Taking into account the conformal invariance of causal properties, this simple result is applicable to some interesting classes of spacetimes which can be written as warped products.

In particular, the {\em standard static} spacetimes are written as $\M=\R\times M, g=-\Lambda d\tau^2+g_0$ where $g_0$ is Riemannian and $\Lambda>0$ depends only on the $M$ part (this is a subclass of the standard stationary spacetimes already introduced, which will  be developed next). Clearly $g$ is conformal to 
a product with $g_R=g_0/\Lambda$. 

Another conformal class are the so-called {\em Generalized Robertson Walker (GRW)} spacetimes \cite{ARS}, $\M = I\times M, g=-dt^2+f^2 g_R$, where $I\subset \R$ is an interval and $f>0$ depends only on the $I$ part. Putting $d\tau=dt/f$, the metric $g$ is conformal to a product spacetime  on $J\times M$, where $J\subset \R$ is another interval computable from $\int dt/f(t)$.  It is easy to check that the causal properties in Proposition \ref{p_escala_causal} remain valid if $\R\times M$ is replaced by $J\times M$, in contrast with the properties of the boundaries (see \S \ref{s5}).
\end{remark}

\subsection{Finsler applications to stationary  spacetimes}\label{s3c}
Next, our aim is to obtain a result which extends Proposition \ref{p_escala_causal} to the standard stationary spacetimes $(\M=\R\times M,g)$ introduced in \eqref{statmet}. 

\subsubsection{Background on the (non-symmetric) Finslerian distance $d_F$} In order to state our main results, let us summarize first some metric properties valid for any Finsler metric $F$ following \cite{CJS11}. Recall from \eqref{dF} that $F$  provides  a (non-necessarily symmetric) distance $d_F$ as well as the corresponding forward and backward balls $B_F^\pm(p,r)$ in \eqref{B+} and \eqref{B-}. These distances yield a natural notion of {\em forward and backward Cauchy sequences} and, then of {\em forward and backward completeness}. Such notions are nicely interpreted taking into account 
 the {\em reversed Finsler metric}, $\tilde F(v) \coloneqq F(-v)$ for all $v\in TM$. Indeed, the latter clearly satisfies
$B_{\tilde F}^+(p,r)=B_F^-(p,r)$. Moreover, if $\gamma$ is a geodesic for $F$, then its reversed parametrization (i.e., $\tilde \gamma(s) \coloneqq \gamma(-s)$ for all $s$) becomes a geodesic for $\tilde F$ and so, one speaks of {\em forward (resp. backward) geodesic completeness} when  all the inextendible geodesics have an upper (resp. lower) unbounded domain. A version of Hopf-Rinow theorem  and some further computations show the following. 

\begin{theorem}\label{hopfrinow} Let $F$ be any Finsler metric  on $M$.

\smallskip

\noindent (A) The following conditions are equivalent (see for example \cite[Theorem 6.6.1]{BCS}):

(A1) forward (resp. backward) completeness of 
$d_F$, 

(A2)  compactness of all its closed forward  (resp. backward) balls, 

(A3)  forward (resp. backward) geodesic completeness  of $F$.

\smallskip


\smallskip

\noindent (B) Let  $d_F^s(p,q) \coloneqq (d_F(p,q)+d_F(q,p))/2$ for all $ p,q\in M$, then $d_F^s$ is a (usual) distance on $M$, called the {\em symmetrized distance}, and the following conditions are equivalent (see \cite[Proposition 2.2, Theorem 5.2]{CJS11}):

(B1) compactness of all the closed balls for the symmetrized distance $d_F^s$, 

(B2)  precompactness of  all the intersections 
$ B_F^+(p,r)\cap  B_F^-(q,r')$, $p, q\in M$, $r, r'>0$, 

(B3) the last property putting $p=q, r=r'$.

 In this case, $d_F^s$ is complete and $F$ (and so $\tilde F$) is convex (i.e., every two points are connected by a minimizing geodesic).
Moreover, the above properties hold if either $d_F$ or $d_{\tilde F}$ is complete.


\end{theorem}
\begin{remark}  (1)  The symmetrized distance $d^s_F$ is not associated with any length space, in general. This underlies the fact that the converse of the last assertion in the  proposition does not hold, that is, there are even Randers examples with $d_F^s$  complete but non-compact closed balls (see \cite[Example 2.3]{CJS11}).

(2) It is easy to check that the forward completeness of $F$ does not imply the backward one. In the particular case of a Fermat metric  \eqref{Fermat}, the forward or backward completeness of $F$ implies the completeness of $-\frac{1}{\Lambda}\tilde g=\frac{1}{\Lambda^2} \omega^2+\frac{1}{\Lambda} g_0$ but the converse does not hold (see \cite[Proposition 5.2 and Example 2.3]{CJS11}).

\end{remark}

\subsubsection{Finslerian description of causality}
It is not difficult to realize that the chronological and causal futures and pasts  in a standard stationary spacetime, $I^\pm(p_0), J^\pm(p_0), p_0=(t_0,x_0)\in \R\times M$,  can be described in terms of the balls for the distance $d_F$ associated with the Fermat metric $F$ defined in \eqref{Fermat} as follows (see \cite[Proposition 4.2]{CJS11}):
$$
I^+(t_0,x_0)= \cup_{s>0} \{(t_0+s)\times B_F^+(x_0,r)\}.
$$

\begin{remark}
This can be seen as a consequence of a further interpretation of the Zermelo problem described in \S \ref{s2}  (see especially Remark \ref{lightZermelo}).  We can consider that the prescribed velocities are the top speeds that a moving object can reach. Accordingly, the trajectories in the spacetime will be causal curves, being lightlike when a certain path is traveled at the top speed. Moreover, by the causal properties, a curve is a solution of Zermelo problem if and only if it remains in $J^+(t_0,x_0)\setminus I^+(t_0,x_0)$, which implies that it is a lightlike pregeodesic and its projection is a minimizing pregeodesic of $F$  of length $t_1-t_0$  (recall Proposition \ref{lightlikegeo}). As in Riemannian Geometry, the geodesics of $F$ are always locally minimizing and the lightlike geodesics departing from $p$ lie initially in $J^+(p)\setminus I^+(p)$.
\end{remark}
These are the keys for the description of  causality using  the Fermat metric $F$ in \eqref{Fermat}, allowing one to prove the following result (see \cite[Thms. 4.3 and 4.4]{CJS11}).
\begin{theorem}\label{t_escala_causal} For any standard stationary spacetime $(\M,g)$ with associated Fermat metric $F$ (as in \eqref{statmet}, \eqref{Fermat} above), the canonical projection $t: \R\times M \rightarrow \R$ is a temporal function and

\smallskip

\noi (A) The following properties 
are equivalent:

(A1) $(\M, g)$ is causally simple.

(A2) $F$ is convex (i.e., each $p,q\in M$ can be joined by a minimizing $F$-geodesic).

\smallskip

\noi (B) The following properties 
are equivalent:

(B1) $(\M, g)$ is globally hyperbolic.

(B2) the intersections 
$ B_F^+(p,r)\cap  B_F^-(q,r'), p, q\in M, r, r'>0$ are  precompact. 

\smallskip

\noi (C) The following properties 
are equivalent:

(C1) The slices $\{t_0\}\times M$ are Cauchy hypersurfaces (so, $t$ is Cauchy temporal).

(C2) The Fermat metric $F$ is forward and backward complete.

\end{theorem}
The comparison between Proposition \ref{p_escala_causal} and Theorem \ref{t_escala_causal} shows that the case of standard stationary spacetimes is much subtler than the case of product spacetimes. Much of this subtlety comes from the following fact. Both cases are presented with a natural slicing by spacelike hypersurfaces which are invariant by the timelike Killing vector field $\partial_t$. In the product case, this slicing is priviledged,  as it corresponds with the integral hypersurfaces of the  distribution $\partial_t^\perp$, orthogonal to $\partial_t$.  In the standard stationary spacetime, however, the slices are not especially priviledged, and one can always  express  the spacetime as a standard stationary one in many different ways, as emphasized next. 

Any function $f: M \rightarrow \R$ which yields a spacelike graph $S^f=\{(f(x),x), x\in M\}\subset \M$  generates a new splitting of the spacetime $(\M,g)$ as a standard stationary one such that $\M\cong \R \times S^f$ (just moving $S^f$ with the flow of $\partial_t$). The Fermat metric associated with this new splitting is $F^f \coloneqq F-df$ (see \cite[Proposition 5.9]{CJS11}). 
 The characterization of the causality conditions in Theorem \ref{t_escala_causal} together with the last observation about the relation between Fermat metrics associated with different slices of the standard stationary spacetime has some striking consequences for Randers metrics.
\begin{cor}\label{appRanders}
Let $(M,F)$ be a Randers manifold (recall \eqref{randers}) and assume that the intersections $B^+_F(p,r)\cap B^-(p,r)$ are precompact for all $p\in M$ and $r>0$. Then  $F$ is convex and  there exists a function $f:M\rightarrow \R$ such that $F^f=F-df$ is a complete Randers metric on $M$.
\end{cor}
Both properties follow from well-known causal properties. Indeed, as stated in Theorem \ref{t_escala_causal}, the precompactness of the intersection of the balls of the Fermat metric is equivalent to the global hyperbolicity of the spacetime, which implies the causal simplicity, and then again by Theorem \ref{t_escala_causal}, the convexity of the Fermat metric. Moreover, as stated in Theorem \ref{t_Cauchytemporal}, if $(\R\times M,g)$ is globally hyperbolic, then it admits a smooth spacelike Cauchy hypersurface, and the Fermat metric associated with this slice 
can be expressed as $F^f=F-df$ for a certain $f:M\rightarrow \R$ and it is complete (see \cite[Theorem 5.10]{CJS11}). 
\begin{remark}\label{appFinsler}
The previous corollary holds for any Finsler metric. This can be proved either directly, as it was done by Matveev in \cite{Matveev}, or generalizing Theorem \ref{t_escala_causal} to stationary Finsler spacetimes as in \cite{CapSta16,CapSta18} (and, then, using the generalization of the existence of a smooth Cauchy spacelike hypersurface for cone structures, firstly obtained by Fathi and Siconolfi in \cite{FS12}). 
\end{remark}

\begin{remark}
  Theorem \ref{t_escala_causal} can be used to get some  different types of results for standard stationary spacetimes: 
\begin{enumerate}[(i)]
\item  As a first consequence, one can obtain some multiplicity results for lightlike geodesics (periodic or between two fixed points) of 
globally hyperbolic stationary spacetimes with arbitrary big arrival times. This follows using multiplicity results for Finsler metrics (see \cite{BiJa11,CJM10,CJM11}).  Moreover,  using the auxiliary product spacetime $((\R\times M)\times \R, g+dr^2)$, one can also obtain multiplicity results for timelike geodesics with prescribed proper time between a point and a vertical line of $(\R\times M,g)$. Finally, using some results on convex boundaries for Finsler metrics in \cite{BCGS11}, it follows the existence of lightlike or timelike geodesics with prescribed proper time for stationary spacetimes with suitable boundaries (see \cite{CGS11}). 

\item   Theorem \ref{t_escala_causal} can also be generalized to pre-Randers metrics, namely, Randers metrics as in \eqref{randers} with arbitrary $\omega$, this time with some subtleties as the slices $\{t_0\}\times M$ do not have to be necessarily spacelike (see \cite[\S 3.2]{HJ21}). This more general case has applications to the study of magnetic geodesics and also to strengthen results about the existence of $t$-periodic lightlike geodesics of stationary-complete spacetimes (see \cite[Corollary 5.12]{HJ21}).  

\item Observe that the Finsler metric	$F^f$ in Corollary \ref{appRanders} has the same pregeodesics as $F$. Moreover, their distances are straightly related having the same triangular function $T(x,y,z)=d(x,y)+d(y,z)-d(x,z)$. As a further interplay between Lorentz and Finsler geometries, the maps which preserve this triangular function are called {\em almost isometries} and they can be identified with conformal maps of the spacetime that preserve $\partial_t$ (see \cite[Proposition 4.7]{JLP15} for details). Moreover, the almost isometries have been characterized in terms of the semi-Lipschitz functions of $(M,F)$ (see \cite{CJ17,DJV20}).
\end{enumerate}

\end{remark}

\subsection{Application to general SSTK spacetimes}\label{s3d}
The previous setting can be extended to any SSTK spacetime, as proven in \cite{CJS22}. This is more technical and only the main ingredients will be sketched.  For the basic notions on SSTK spacetimes, recall \S \ref{SSTKsec}. 


\begin{remark}\label{r_killinghorizons}
In Relativity, SSTK spacetimes are applicable to describe stationary black holes.\footnote{See \cite{Wald} for background, especially Chapter 12 and \S 12.3, \S 12.5.} Indeed, the so-called {\em Killing horizons}   can be written  as  degenerate hypersurfaces of an SSTK spacetime determined by $\Lambda=0$  (see Figure \ref{KillingHor}).  Typically,  this notion appears under assumptions on asymptotic flatness (so that one can think of Lorentz-Minkowski observers at infinity) and $\Lambda$ becomes positive  asymptotically. Thus, the open subset $\Lambda> 0$ is expected to have a  connected part which includes the asymptotic region. This would be the region outside the Killing horizon and (each connected component of) its boundary would be a Killing horizon.

 \begin{figure}
 \centering
 \begin{tikzpicture}[scale=0.9]
  \draw (0,0) ellipse (4cm and 2cm) node[text width=3.1cm]{\centering Black hole region \\ $g(K,K)(=-\Lambda)>0$};
  \draw[red] (0,1.6) circle (0.4cm);
  \draw[blue,->,thick] (0,2) -- (0,1.6);
   \draw[red] (0,-1.6) circle (0.4cm);
   \draw[blue,->,thick] (0,-2) -- (0,-1.6);
   \draw[red] (-3.6,0) circle (0.4cm);
   \draw[blue,->,thick] (-4,0) -- (-3.6,0);
   \draw[red] (3.6,0) circle (0.4cm);
   \draw[blue,->,thick] (4,0) -- (3.6,0);
  \end{tikzpicture}
  \caption{Killing horizon determined by the region $\Lambda=0$ in a slice $\{t_0\}\times M$. The blue arrows correspond to the wind $ W $ in the critical region for Zermelo data as in Sect. \ref{subsubsec:wind_Riem}.}
 \label{KillingHor}
 \end{figure}
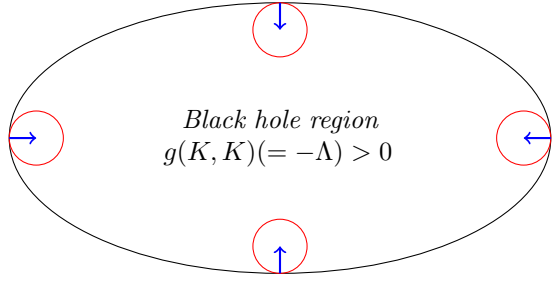

In Kerr spacetime,  the so-called  {\em ergosphere} is its natural  Killing horizon. 
In the case of  Schwarzschild spacetime, such a horizon is equal to its {\em event horizon}, informally, the hypersurface that, once crossed by an observer coming from infinity, cannot be crossed again (a spacetime is expected to be globally hyperbolic outside its event horizon). Indeed, Schwarzschild spacetime admits a Killing vector field $K$ which is timelike  outside this horizon and makes it a static spacetime there (thus, this Killing $K$ becomes priviledged because its orthogonal distribution is integrable; recall Remark \ref{r_static}). The static description  fails in the horizon, but Schwarzschild  spacetime  can be  described as an SSTK therein (being written the region outside the horizon  as a  stationary spacetime).  
\end{remark}

\subsubsection{The particular case $\Lambda\geq 0$}\label{Sranders-kropina}
 Observe that when $\Lambda\geq 0$,  the metric $F$ in \eqref{fermatArbitraryF} is always of Randers type as in \eqref{randers} (when $\Lambda>0$) or of Kropina type  $F(v)=-g_0(v,v)/(2\omega(v))$ (when $\Lambda=0$). Accordingly, we will say that $F$ is a Randers-Kropina metric in this case. 
 A noticeable subtlety is that the domain of $F$ (and, so, the velocities of the admissible curves joining two points) is restricted in the Kropina case to the half space $\omega<0$   (recall \eqref{conicA}).  So, the {\em $F$-separation} $d_F$,  defined formally as a Finslerian distance,   will yield $d_F(x,y)=\infty$ if there is no admissible curve from $x$ to $y$.  Even though not all the statements in Hopf-Rinow theorem hold now, $d_F$ is fairly well-behaved.\footnote{In particular, it is continuous outside the diagonal \cite[Theorem 4.5]{CJS22} (this property does not hold for any conic Finsler metric \cite[Proposition 3.9]{JS_Pisa}).} Indeed, regarding $d_F$ as the $F$-separation one can ensure  (see \cite[Theorem 4.9]{CJS22}): 
\begin{quote} {\em The characterizations of  causality and Cauchy hypersurfaces of a standard stationary spacetime in Theorem \ref{t_escala_causal} remain valid for any SSTK spacetime with $\Lambda\geq 0$.} \end{quote}

\subsubsection{The general case} \label{3.5.2}

 Recall that when $\Lambda<0$ there are two pseudo-Finsler\footnote{ Here pseudo-Finsler means that the fundamental tensor \eqref{fundten} is non-degenerate and the domain is conic, not necessarily the whole $TM$; in particular, so are the Finsler metric $F$   and  the Lorentz-Finsler one $F_l$. } metrics $F$ and $F_l$ (see \eqref{fermatArbitraryF} and \eqref{fermatArbitraryFl}) both defined in the conic subset
\[A_l=\{ v\in TM: \Lambda<0, \ \omega(v)^2+\Lambda g_0(v,v)>0, \ -\omega(v)>0\},\]
and all the lightlike geodesics of  $(\R\times M,g)$ project into pregeodesics of $F$ or $F_l$, or into lightlike pregeodesics of $\omega^2+\Lambda  g_0$ 
(which matches \eqref{-tildeh} up to a sign; see Theorem \ref{lightlikegeo2}).   Then 
a precise description of the causality and Cauchy hypersurfaces of SSTK spacetimes in terms of its associated wind Riemannian structure is available; see \cite[Theorem 5.9]{CJS22}. 
Even though it is subtler than the previous cases,\footnote{Among the subtleties appearing here, it is worth mentioning the step {\em causal continuity} (i.e., the spacetime is  causal and the causal future and past  vary continuously with the point) in the ladder of spacetimes. This  is more restrictive than stable causality and it is satisfied by all the standard stationary ones, as well as whenever $\Lambda\geq 0$; however, it is not satisfied by all the SSTK spacetimes.} to determine when the slices $t=$ constant  are Cauchy hypersurfaces  is simple. Indeed, only the conic metric $F$ (which is defined on the whole spacetime) becomes relevant. Therefore, as emphasized in \cite{JS_proc}, the  extension $\bar F$ of $F$ to all the projections of future-directed causal vectors is continuous and permits an extension 
$d_{\bar F}$ of the $F$-separation  on the whole $M$ working as in \S \ref{Sranders-kropina}
(see  \cite[\S 3.1.1 and \S 3.2.1]{JS_proc}). Then the  Cauchy completeness of $d_{\bar F}$ becomes naturally equivalent to the completeness of the wind Riemannian structure of the SSTK spacetime \cite[Theorem 3.23]{JS_proc} and, as a consequence \cite[Corollary 4.1]{JS_proc}:
\begin{quote}
{\em For any SSTK spacetime, the slices $t=$ constant are Cauchy if and only if the extended distance $d_{\bar F}$  is (forward and backward) complete.}
\end{quote}

\section{Finsler applications}\label{s4}

Next, some applications of the Lorentz-Finsler approach to the Finslerian setting will be gathered. We have already pointed out some of them. Indeed,   one of the first consequences  was to highlight the importance of the compactness of the symmetrized closed balls, or equivalently, the precompactness of the intersection of forward and backward balls. This condition is equivalent to the global hyperbolicity of the stationary spacetime associated with a Randers metric. Moreover, it also  implies the geodesic connectedness of $(M,F)$ (recall Theorem \ref{hopfrinow})  and the existence of a complete Finsler metric $F^f=F-df$ for a certain function $f:M\rightarrow \R$ with the same pregeodesics as $F$ (recall Corollary \ref{appRanders} and Remark \ref{appFinsler}). 

The subtlety of this property was stressed in \cite[Example 2.11]{JS_proc} by exhibiting a Randers manifold $(M,R)$ satisfying:  {\em $(M,R)$ does {\em not} have compact closed symmetrized balls but its universal covering $(\hat M,\hat R)$ does}. This follows by applying Theorem \ref{t_escala_causal} to a globally hyperbolic stationary  spacetime  previously  constructed by Harris \cite[Example 3.4(b)]{Harris15}, which admitted a non-globally hyperbolic quotient by isometries preserving the timelike Killing field. 

A further application of stationary spacetimes to Randers metrics is a result about the smoothness of the distance to (or from) a closed subset. As it was shown in \cite[Theorem 5.12 and Corollary 5.13]{CJS11}, the subset of points where this distance is not smooth must have negligible $n$-Hausdorff measure.


\subsection{Randers-Kropina metrics}  The applications of spacetimes to Finsler Geometry can be even more powerful when one considers Randers-Kropina metrics or, more generally, wind Riemannian structures. This is because these metrics have certain type of singularities in the boundary of the domain, which disappear from the spacetime viewpoint. First, we have a direct extension to some properties for Finsler metrics stated above (see \cite[Proposition 6.6 and Corollary 6.11]{CJS22}):
\begin{theorem}
Let $(M,F)$ be a Randers-Kropina manifold (recall \S \ref{Sranders-kropina}) such that all the intersections of its backward and forward balls are precompact. Then
\begin{enumerate}[(i)]
\item If $d_F(p,q)<+\infty$, (or equivalently, there exists an admissible curve between $p$ and $q$), then there is a minimizing geodesic joining them.
\item There exists a smooth function $f:M\rightarrow \R$ such that $F-df$ is a geodesically complete Randers-Kropina metric.
\end{enumerate}
\end{theorem}
 For arbitrary wind Riemannian structures, it is possible to give a version of part $(i)$ of the last theorem with a more technical definition of forward and backward balls (see \cite[Definition 2.26]{CJS22} and \cite[Proposition 6.6]{CJS22}). 
The extension of part (ii) is not straightforward. Again, the precompacteness of the intersections of the forward and backward balls imply global hyperbolicity and, thus, the existence of a Cauchy hypersurface $S$. However, as the integral curves of the Killing field $K=\partial_t$ may be spacelike, they can cross $S$ more than once or not at all (thus, an SSTK splitting on $S$ would not be obtained).

In \cite[Theorem 1.5 (B)]{CMMM} the authors prove geodesic connectivity  for the Kropina manifold $(M,g/\omega)$ with $M$ compact and $\omega\wedge d\omega\not=0$. This follows from part $(i)$ above by observing that $\omega\wedge d\omega\not=0$ implies that every two points are connected by an admissible curve. Moreover, stationary spacetimes have also been used to obtain multiplicity results for Kropina geodesics between two points in a compact manifold (see \cite{CGMS}).

\subsection{Classification of Randers and wind Riemannian spaceforms} 
Let us recall that flag curvature is a fundamental invariant of a Finsler manifold which measures how geodesics with initial velocity in a certain plane spread apart. Flag curvature is the generalization of sectional curvature to Finsler Geometry and the problem of classifying the constant flag curvature Finsler manifolds remains as one of the biggest challenges. This classification has only been achieved for some families of Finsler metrics such as Randers and Kropina \cite{BaRoSh04,YosSab12}.
Let us recall that Zermelo problem was essential to obtain the final classification of Randers spaceforms in the celebrated paper \cite{BaRoSh04}. Indeed, it turns out that a Randers manifold has constant flag curvature if and only if its Zermelo data $(\ho,W)$ is given by a Riemannian metric $\ho$ with constant sectional curvature and a homothety $W$ of $\ho$. Differently from the Riemannian counterpart, many of these Randers spaceforms are not geodesically complete. This anomaly disappears when one   drops the Randers restriction $g_R(W,W)<1$ and, thus, considers the  wider family of wind Riemannian structures. Recall that there is a natural notion of geodesic intrinsic to such a structure $\Sigma$ \cite[Definition 2.35]{CJS22} and, then, of geodesic completeness, the latter implying the geodesic completeness of both its conic Finsler metric $F$ and the Lorentz Finsler one $F_l$ \cite[Corollary 5.6 and Figure 10]{CJS22}. Moreover,  as the fundamental tensors of $F$ and $F_l$ are computable and non-degenerate away from the boundary of $A_l$, their flag curvatures are well-defined; 
 so, $\Sigma$ is called of constant flag curvature when the flag curvatures of both  $F$ and $F_l$ are equal to the same constant.  Then, the techniques used to prove the Randers classification can be extended to these structures to obtain the following (see \cite[Theorem 3.12]{JS17}).   
\begin{theorem}\label{e_cfc}
The complete simply connected wind Riemannian structures with constant flag curvature lie in one of the following two exclusive cases, determined by the Zermelo data $(g_R,W)$:

(i) $(M, g_R)$ is a model space of constant curvature and $W$ is any of its Killing
vector fields.

(ii) $(M, g_R)$ is isometric to $\R^n$ and $W$ is a properly homothetic (non-Killing)
vector field.

Moreover, the inextensible simply connected Randers (resp. Kropina) manifolds with constant flag curvature are the maximal simply connected open subsets of the previous wind Riemannian structures where the wind is mild (resp. critical).   
\end{theorem}
Notice how, in the last assertion, the assumption of completeness must be replaced by the inextensibility of the Randers or Kropina metric as a manifold of this same type.
We can go further, and give a characterization in terms only of the conic metric $F$ defined on the whole  $M$.  Indeed, it is natural to rename it as the {\em Zermelo metric} $Z$ for the data $(g_R,W)$ (given explicitly  in \eqref{zermeloArbitrary}). 
 At the end of \S \ref{3.5.2}, we saw that such a  $Z$ determines univocally an extended distance  $d_{\bar Z}$. Moreover, the Cauchy completeness of the latter is equivalent to the  completeness of the wind Riemannian structure (see \cite[Theorem 3.23]{JS_proc}), thus yielding the following. 
\begin{cor}
A   Zermelo metric $Z$ with complete extended distance $d_{\bar Z}$  on a simple connected manifold $M$ has constant flag curvature $k\in\R$  if and only if its Zermelo data  $(g_R,W)$ satisfies either (i) or (ii) in Theorem \ref{e_cfc}.
\end{cor}

\section{Interplay Finsler/ Lorentz for boundaries} \label{s5}

Let  $(M,g_R$) be a Riemannian manifold. Its elementary Cauchy boundary $\partial_C M$ provides a completion $M_C$ and, when $g_R$ is complete, its Gromov boundary $\partial_G M$  (see~\cite{Gr}) provides a compactification $M_G$. 
In the case of Hadamard manifolds, this compactification agrees with the previous one by  Eberlein and O'Neill (see~\cite{EO}), which was introduced in a very different way by using Busemann functions associated with rays.\footnote{Hadamard manifolds (which are complete, simply connected and with non-positive sectional curvature) become diffeomorphic to $\R^n$ by Cartan-Hadamard theoRemark Eberlein and O'Neill's boundary becomes a topological sphere $S^{n-1}$ located  at infinity. Indeed, the Busemann functions yield a quotient in the space of rays so that  each function can be regarded as a direction at infinity; then, the set of all these directions can be regarded as the sphere  $S^{n-1}$.} 
Being the main properties of these boundaries well established since the beginning of the eighties,  natural questions about the relation among them, as well as its  extension to (possibly non-reversible) Finslerian metrics, had remained dormant. However, an additional motivation for their study came from the links with the Lorentzian setting. 

Roughly, the link with the Riemannian case appears when one computes the causal boundary of a product $\R\times M, g=-dt^2+g_R$. In a natural way, the computation of this boundary leads to consider a sequential compactification $M_B$ of $(M,g_R)$ in terms of Busemann-type functions, thus extending Eberlein and O'Neill compactification.  This compactification includes the Cauchy completion $M_C$ in a natural but subtle sense, and it is also related to  Gromov's  $M_G$. 

Amazingly, 
$M_B$ and $M_G$ are equivalent except in some rather pathological cases, which also correspond with known pathologies of the causal boundary. 
All this can be extended to the Finslerian setting by providing a natural link with the causal boundary of  standard stationary spacetimes (eventually extendible to Finslerian spacetimes). 
However, the possible non-reversibility of 
the Finslerian metric introduces subtleties 
even at the level of the Cauchy boundary.
The systematic analysis 
of all these issues 
was  carried out in \cite{FHS_Memo}, to be followed next.

\subsection{Gromov compactification for incomplete  Finslerian manifolds} 
Let us start 
considering a metric space $(M,d)$ associated with a connected Riemannian manifold $(M,g_R)$ or, with more generality, $M$ can satisfy just to be connected, locally compact 
and second countable,  while the distance $d$ is just derived from a length space.\footnote{Notice that, for a smooth curve $c$  in such a space,  the role  of its Riemannian norm $g_R(\dot c(s),\dot c(s))$ at each $s$ is played by the local  dilatation there; moreover, notions such as geodesic or cut point have a natural  sense (see \cite[Chapter 1]{Gr99} for background).} In particular, 
the metric space associated with any {\em reversible} Finsler manifold is included now and, later, we will refer  to the Finsler case when  taking into account non-reversibility. 

\subsubsection{The symmetric $d$ case} The Cauchy completion $M_C$ of such a space is standard, and we emphasize  that $M_C$ may be non-locally compact (thus, it will not lie under the general hypothesis for $(M,d)$ above). Indeed, it is not difficult to construct a  bidimensional Riemannian example starting at a variation of the comb space
\begin{equation}
\label{e_comb1}  M=\left( (0,\infty)\times \{0\} \right) \cup \left( K\times (0,1)\right) \subset \R^2, \quad  \hbox{where} \quad K=\{1/m: m\in \N \}.
\end{equation}
Here,  $(0,0)$  is identifiable to a point of $\partial_CM$ which does not have any compact neighborhood in $M_C$.

Let us construct the Gromov compactification of $(M,d)$ without the usual assumption on completeness for $d$. Consider the space  of all the 1-Lipschitz functions $\mathcal{L}_1(M,d)$. For each $x\in M$, the function $M\ni y\mapsto d_x(y) \coloneqq d(x,y)$ belongs to $\mathcal{L}_1(M,d)$. Moreover,  $d_x$ as well as any function $d_x+C$, where $C$ is a constant, determines univocally $x$. Thus, $M$ can be identified with a subset of the quotient space $\mathcal{L}_1(M,d)/\R$ under the relation of equivalence:
\begin{equation}
\label{e_relac_equiv}
f\sim f' \Leftrightarrow f-f'=C \in \R, \quad \hbox{where} \quad f, f'\in  \mathcal{L}_1(M,d).
\end{equation}

\begin{defi} The Gromov completion  $M_G$ of $(M,d)$ is the closure of $M$ in $\mathcal{L}_1(M,g)/\R$ (with, say, the  uniform convergence on compact sets topology).
\end{defi}
It is not difficult to check that $\mathcal{L}_1(M,g)/\R$ is compact and, moreover  (see \cite[Theorem 4.12 and Corollary 4.13]{FHS_Memo}): 
\begin{prop}\label{p_Gromov}
	$M_G$ is a compact
	metrizable space and $M_C$ is naturally included in it.
	The inclusion $M_C\hookrightarrow M_G$ is continuous, and it is an embedding if and only if $M_C$ is locally compact; moreover,   $M\hookrightarrow M_G$ is  a dense embedding. 
\end{prop}
To check that the inclusion of $M_C$ in the example \eqref{e_comb1} is not an embedding, notice that the sequence $\{(1/n,1/2)\}_n$  converges to $(0,0)$ in $M_G$ but it is not convergent in $M_C$. This suggests some subtleties for this boundary. Indeed,  Gromov's boundary $\partial_G M \coloneqq M_G\setminus M$ 
is divided into a {\em Cauchy-Gromov} boundary 
$\partial_{CG}M$, whose points are the limits of bounded sequences in $M$, and a {\em proper} Gromov boundary $\partial_{\mathcal G}M$ containing limits of unbounded sequences.  Clearly, $\partial_{CG}M$ contains the Cauchy boundary $\partial_C M$; however it may contain more points if $M_C$ is not locally compact. 

\begin{exe}\label{e_pathological}  Modify  the space in \eqref{e_comb1} by adding an upper half line:
$$	M=\left( (0,\infty)\times \{0,1\} \right) \cup \left( K\times (0,1)\right) \subset \R^2, \quad  \hbox{where} \quad K=\{1/m: m\in \N \}.
	$$
	Now, for each $y\in (0,1)$ the  sequence $\{(1/n,y)\}_n$
	converges to a distinct limit, so that  $\{0\}\times [0,1]$ can be regarded as a subset of $\partial_{CG}M$. We emphasize that {\em no boundary point $(0,y), y\in(0,1)$, can be an endpoint of a curve starting at $M$}.
\end{exe}

\subsubsection{The non-symmetric $d$ case} When considering the Finsler case, $d$ will mean the non-necessarily symmetric distance defined in \eqref{dF}, and the following subtleties 
must be taken into account. Following   
\cite{Z}, a map $d: M\times M \rightarrow \R$ is called a {\em generalized distance} when it satisfies, for all $x, y, z \in M$: (a1)~$d(x, y) \geq 0$, (a2)~$d(x, y) = d(y, x) = 0$ if and only if $x = y$, (a3)~$d(x, z) \leq  d(x, y) + d(y, z)$, and
(a4)~given a sequence $\{x_m\}_m \subset  M$ and $x \in M$, then: 
$\lim_{m\rightarrow \infty} d(x_m, x) = 0$ if and
only if $\lim_{m\rightarrow \infty} d(x,x_m) = 0$; when the hypothesis (a4) is dropped, then $d$ is called a {\em quasidistance}. 
A generalized distance $d$ gives two notions of Cauchy sequences (forward and backward) and, then, forward and  backward Cauchy completions $M_C^+, M_C^-$, respectively.\footnote{There is a non-equivalent way to define forward (and backward) Cauchy sequences; however, it would yield the same forward Cauchy completion (see \cite[Sect. 3.2.2]{FHS_Memo}).} Moreover,  the distance $d^s$ obtained by symmetrizing $d$ provides another completion $M^s$ and the corresponding boundaries satisfy $\partial^s_CM=\partial^+_CM\cap \partial^-_CM$ in a natural way.   The continuous extension of $d$ to such a completion is  only a quasidistance; indeed, the topologies generated by the forward and backward balls are not equivalent and $M_C^+$ is only a $T_0$ space. Notice that, in our previous study of  Finsler metrics, we  considered a non-symmetric distance  which is, indeed, a generalized distance, and all the assertions above apply \cite[Chapter 3]{FHS_Memo}. 

In order to consider the Gromov completion for a Finsler manifold, notice that there are two non-symmetric notions of 1-Lipschitzian:

$$ \begin{array}{rl}
\mathcal{L}_1^+(M,d)= & \{f \, \hbox{smooth}: \, f(y)-f(x)\leq
d(x,y)\}, \\
\mathcal{L}_1^-(M,d)= & \{f \, \hbox{smooth}: \, f(x)-f(y)\leq d(x,y)\}.
\end{array}$$
Accordingly, there are two Gromov compactifications $M_G^\pm$,\footnote{Its consistency relies on a non-symmetric version of Arzela theorem \cite[Theorem 5.12]{FHS_Memo}.}
namely:  
$M_G^+$ is the closure of $M$ in $\mathcal{L}_1^+(M,d)/\R$.  

\begin{remark} The inclusion $M_C^+ \hookrightarrow
	M_G^+$ is subtler than in Proposition \ref{p_Gromov},  as it  satisfies now (see \cite[Corollary 5.25]{FHS_Memo}):
	
	(A) it is continuous  if and only if  the backward balls generate a
	finer topology on $M_C^+$ than the forward balls,
	
	(B) it is an embedding when: (B1) $M^+_C$ is
	locally compact (as in the Riemannian case) and  (B2) the extension of $d$ to $M_C^+$ is a
	generalized distance.
\end{remark} 

\subsection{The causal boundary of a spacetime}
For spacetimes, there are two conformally invariant boundaries which are applied to general relevant classes of spacetimes. The first one is the so-called {\em conformal boundary}, introduced by Penrose in \cite{P}, which  underlies  notions such as  {\em asymptotic flatness}, and it is widely used in Relativity. Essentially, the idea is to find a suitable open conformal embedding of the spacetime in a bigger one and, then, to regard its topological boundary  as the conformal one. The second one is the {\em causal boundary}, firstly introduced by Geroch, Kronheimer and Penrose \cite{GKP}, but later redefined several times (we will refer  to the last one, \cite{FHS11}). This boundary $\partial_c\M$ is defined in an intrinsic way for any spacetime under the weak condition  of being {\em strongly causal}.\footnote{Intuitively, it does not admit ``almost closed'' causal curves. A formalization of this property is that each point $p\in M$ has a neighborhood $U$ such that any inextendible causal curve starting at $U$ will leave  $U$ at some point so that it will not return to $U$.}
There are general conditions which ensure that these boundaries agree (as well as  counterexamples otherwise); see \cite[Sect. 4 and Appendix]{FHS11}.\footnote{Noticeably, they agree in the class of globally hyperbolic spacetimes-with-timelike-boundary, 
	see~\cite{AFS}.}

\begin{figure}
	\centering
	
	\includegraphics[width=0.45\paperwidth,height=0.16\paperheight]{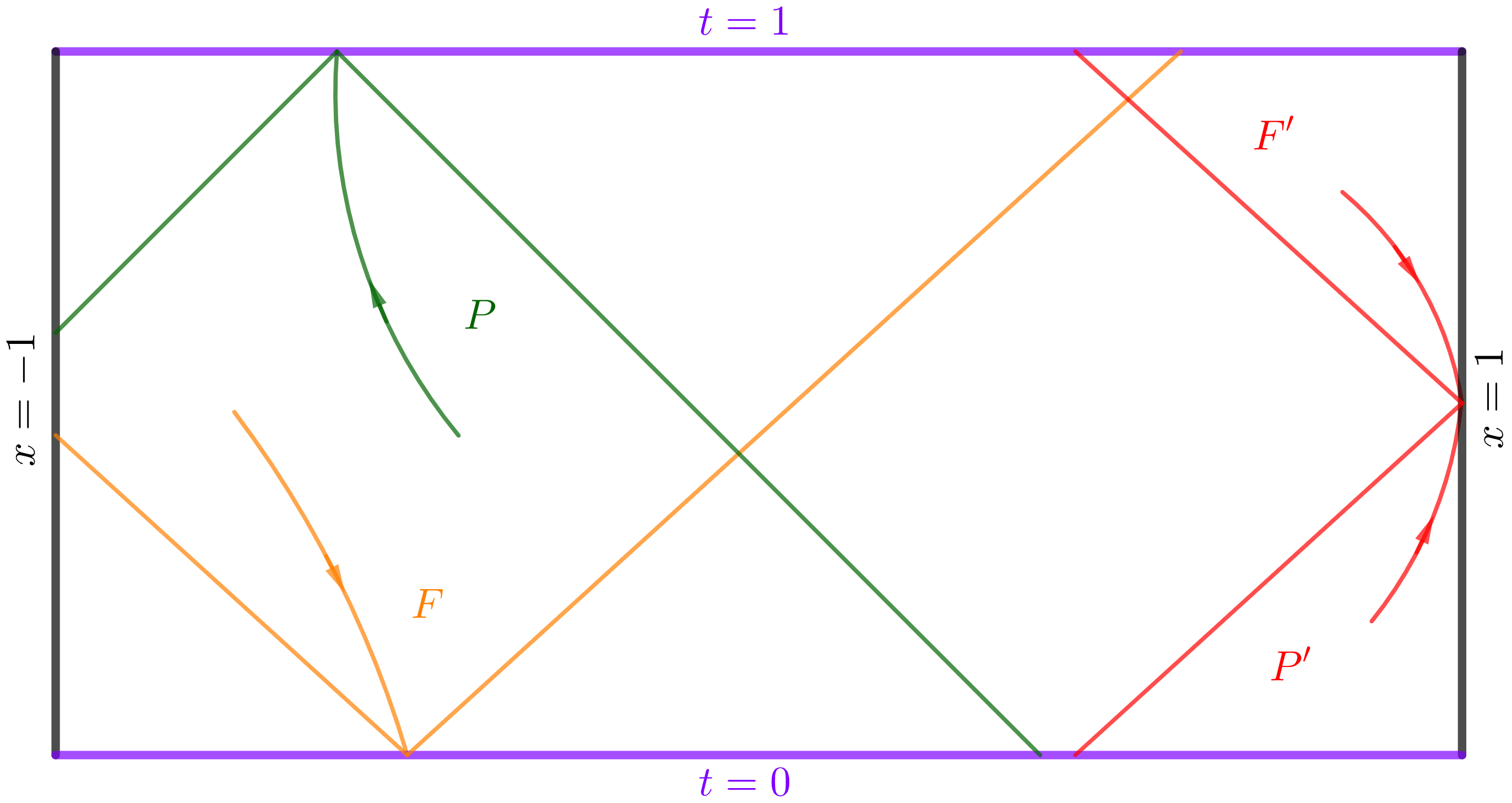}
	\caption{Let $\M= (0,1)\times (-1,1)$ (in Lorentz-Minkowski spacetime). Green $P$ (resp. orange $F$) corresponds with  a TIP (resp. TIF) which represents a point in the causal boundary; intuitively, this point is identifiable with a boundary point in the line $t=1$ (resp. $t=0$). Red $P', F'$ are a TIP and a TIF that intuitively represent the same boundary point at $x=1$. They are paired by Szabados relation, yielding a single point of $\partial_c\M$.}
	\label{Fig0}
\end{figure}

\begin{figure}
	\centering
	
	\includegraphics[width=0.45\paperwidth,height=0.16\paperheight]{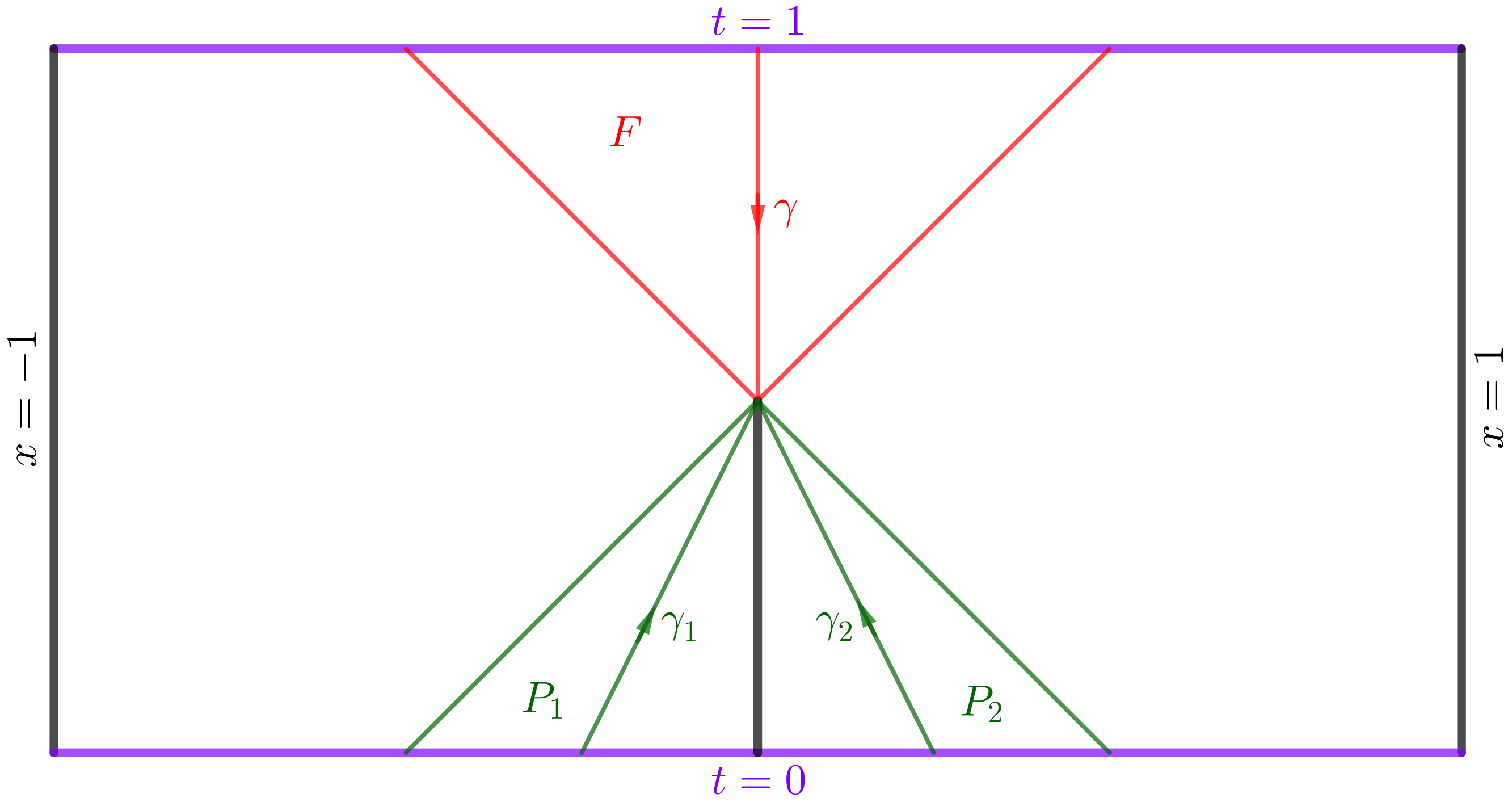}
	\caption{Modify the example in Figure \ref{Fig0} by removing a segment in the $t$ axis, $\M ' \coloneqq \M \setminus \{(t,0): t\leq 1/2\}$. The TIF $F=I^+(\gamma)$, with $\gamma: (0,1/2)\ni t\mapsto (1-t,0)$ is Szabados related with each one of the TIP's $P_1, P_2$, where $P_i=
		I^-(\gamma_i)$, with $\gamma_i: (0,1/2)\ni t\mapsto (t,(-1)^{i} \, ((1/4)-t/2))$. Then,  $P_1\sim F$ and $P_2\sim F$. Therefore, $(P_1, F)$ and $(P_2, F)$ are two distinct points of $\partial_c \M$ (they are not Hausdorff separated by the chronological topology).}
	\label{Fig0b}   
\end{figure}

Here, we are interested in the causal boundary $\partial_c\M$, which will be described  very briefly now, and we refer to \cite{FHS11} for exhaustive details and references. To construct it, one starts defining a {\em terminal indecomposable past} (resp. future) set $P$, or {\em TIP} (resp. $F$, {\em TIF}) for short, as the chronological past (resp. future) of any inextendible future-directed (resp. past-directed) timelike curve $\gamma$, i.e., $P=I^-(\gamma )$ (resp. $F=I^+(\gamma)$). The set of all the TIP's (resp TIF's) is the {\em future} (resp. {\em past}) causal boundary $\partial_c^+\M$ (resp. $\partial_c^-\M$). To construct $\partial_c\M$ one has to take into account that a TIP and a TIF might  represent intuitively the same boundary point (see Figure \ref{Fig0}). So, one introduces the {\em Szabados relation}: $P\sim F$ iff, on the one hand,
$P$ is included in the  common past of $F$ (i.e., ($\cap_{x\in F} I^-(x) \supset P$) and $P$ is maximal among the TIP's satisfying this property, and on the other hand, 
the dual property holds for $F$.
Then $\partial_c \M$ is composed of $\partial_c^+\M \cup \partial_c^-\M$ up to the pairings introduced by $\sim$. Even though this is a neat definition,  examples such as  Figure \ref{Fig0b} show that the pairings may be 
non-trivial. In what follows, all the elements of $\partial_c\M$ will be regarded as a pair $(P,F)$ with the convention that $F=\emptyset$ (resp. $P=\emptyset$) when $P$ (resp. $F$) is unpaired.
\begin{remark}\label{r_bordecausal}
	(1) Globally hyperbolic spacetimes are characterized as the strongly causal ones whose $\partial_c \M$ is composed only of unrelated pairs $(P,\emptyset), (\emptyset, F)$, 
	see \cite[Theorem 3.29]{FHS11}. 
	(2) The chronological relation $\ll$ in $\M$  (introduced around \eqref{e_causalrelations}) admits a natural extension  $\overline{\ll}$ to $\partial_c M$, namely,  $(P,F) \overline{\ll} (P',F')$ whenever $F\cap P'\neq \emptyset$. 
\end{remark} 

The most natural topology for $\partial_cM$ is the so-called {\em chronological topology}. We will not go into the details of this topology, but just point out the following two important features,
in relation to the  general Busemann completion $M_B$ of any Finslerian manifold $M$ to be described below:

\ben \item The topology of $M_B$ is inspired by the chronological topology of the causal completion of a standard stationary spacetime $\M=\R\times M$. In particular,   $\partial_c \M$ will be described directly from the Busemann boundary $\partial_B M$.
\item When this topology is Hausdorff, $M_B$ is identifiable to $M_G$. Otherwise, the non-Hausdorff property of $M_B$ will be related to the appearance of somewhat pathological properties of $M_G$, as the one emphasized in Example \ref{e_pathological} (see Theorem \ref{t_B} below).
\een
This second item supports the relevance of the Busemann completion even in a purely Riemannian setting, and the first one supports  the previously defined causal boundary.


\subsection{A new Busemann boundary}

Next, we will construct a general Busemann boundary which generalizes Eberlein and O'Neill's, following 
\cite[\S 4.2 and \S 5.2]{FHS_Memo}. We will start at the symmetric case with a length metric space $(M,d)$ as in the case of Gromov's,  but  one can consider a connected  Riemannian manifold $(M,g_R)$ to be more specific.

In the standard approach, one assumes  the completeness of $d$, considers a ray $c$ (a half unit geodesic  with no cut locus) and defines its Busemann function $b_c$ as
$b_{c}(x_{0}) \coloneqq \lim_{t\rightarrow \infty}(t-d(x_{0},c(t)))$
for all $x_0\in M$. However, we will drop completeness and admit (generalized) Busemann functions for more general curves, namely: 
\begin{equation}
\label{e_bc}
b_{c}(x_{0}) = \lim_{t\rightarrow \Omega}(t-d(x_{0},c(t)))
\quad \hbox{for any} \; c:[0,\Omega)\rightarrow M \; \text{with} \;  F(\dot c) \leq 1.
\end{equation}
Easily, 
if $b_c$ is $\infty$ at some $x_0\in M$ then $b_c\equiv \infty$, and $B(M)$ will denote the set of all the finite Busemann functions on $M$. 
$B(M)$ will be regarded as a topological space with the 
{\em chronological topology} defined by means of a  limit operator
$L$. Specifically, 
given $\{f_m\}_m\subset B(M)$, the subset $L(\{f_m\}_m) \subset B(M)$
is defined by
$$f\in L(\{f_m\}_m) \Leftrightarrow
\left\{\begin{array}{l}(a)\; f\leq \liminf_{m}f_{m} \quad \hbox{and}\\
(b)\; \forall g\in B(M) \hbox{ with } f\leq g\leq \limsup_{m}
f_{m}, \hbox{ it is } g=f.
\end{array}\right. 
$$
Then, the topology is defined by declaring that a subset $C \subset B(M)$ is  {\em closed} if and only if $L(\sigma ) \subset C$
for any sequence $\sigma= \{f_m\}_m \subset C$.

\begin{defi}
	As a pointset, the  Busemann completion $M_B$ of $(M,d)$ is 
	$M_B \coloneqq B(M)/\R$  (this is the quotient by an additive constant as in \eqref{e_relac_equiv}), endowed with the quotient of the chronological topology on $B(M)$.
\end{defi}
Notice  that $B(M)\subset \mathcal{L}_1(M,g)$ and, thus, naturally $M_B\subset M_G$.
However, the inclusion may be non-continuous (indeed, the topology of $M_B$ is always coarser  than that of $M_G$). 
Observe also that one can consider $M_C\subset M_B$, as the points in the Cauchy completion correspond to  two  Busemann functions for curves with  
$\Omega<\infty$. Indeed, one has the disjoint union 
$M_B=M \cup \partial_C M \cup
\partial_{\mathcal{ B}} M $, where $\partial_{\mathcal{ B}} M$ contains the classes of Busemann functions with $\Omega=\infty$. 
\begin{theorem}\label{t_B} {\em \cite[\S 5.2]{FHS_Memo}}  The following properties hold:
	
	\noindent (A) $M_B$ is sequentially compact.
	
	\noindent (B) $M_B$ is a $T_1$ topological space (and points in the boundary 
	may be non-$T_2$
	related).

	\noindent (C)  $M_C\hookrightarrow M_B \hookrightarrow M_G$ but the topology of $M_B$ is coarser than the others.
	
	
	
	\noindent (D) They are equivalent:
	
	(D1) $M_B=M_G$ as pointsets. 
	
	(D2) $M_B$ is Hausdorff.
	
	(D3) No sequence contained in $M$ converges to two distinct points in $\partial_B M$.

	(D4) $M_B \hookrightarrow M_G$ is an embedding.
	
	(D5) $M_B \hookrightarrow M_G$ is a homeomorphism.
	
\end{theorem}
Recall that, from (D), $M_G$ contains ``extra points'' when $M_B$ and $M_G$ are not naturally equivalent (recall Example \ref{e_pathological} and see also \cite[Theorem 5.39 and Remark 5.41]{FHS_Memo}).

In the Finslerian case,  take into account that, now, each curve $c$ yields two Busemann functions, depending on the ordering of the arguments  of $d$ in \eqref{e_bc}. So, as in the case of Gromov's, we have  two Busemann completions $M_B^\pm$ depending on that ordering. Then, as in the Riemannian case, one has canonical inclusions $M_C^\pm \hookrightarrow M_B^\pm \hookrightarrow M_G^\pm$ 
and Theorem \ref{t_B} is extended naturally to this case,  \cite[Theorem 5.39]{FHS_Memo}. 
In particular, the topology of $M^+_B$ is coarser than those of $M_C^+$ and $M^+_G$, $M^+_B$ is identifiable to $M^+_G$ if and only if $M^+_B$ is Hausdorff, and
the Hausdorffness of $M_B^+$ is independent of that of $M_B^-$. Moreover, one also has
 the disjoint unions:
$$M_B^+=M \cup \partial_C^+ M \cup
\partial_{\mathcal{ B}}^+ M , \quad M_B^-=M \cup \partial_C^- M \cup
\partial_{\mathcal{ B}}^- M. $$

\subsection{The causal boundary of stationary spacetimes } 

Let us describe $\partial_c\M$ for a standard stationary spacetime $\M= \R\times M$ as in \eqref{statmet} by using the Busemann completions $M_B^\pm$ of the Finsler manifold $(M,F)$, where $F$ is the Fermat metric in \eqref{Fermat}. The generalized distance $d$ associated with $F$ will be denoted here $d^+$ (as it will be related to the future causal boundary $\partial_c\M^+$) and its reversed one $d^-$.  Observe  that  the interpretation of the Fermat metric 
provides the following characterization of the chronological relation $\ll$ (recall \eqref{e_causalrelations}):
$$(t_{0},x_{0})\ll (t_{1},x_{1}) \Leftrightarrow
d^+(x_{0},x_{1}) < t_{1}-t_{0}.$$

Let us see how Busemann functions appear when one computes the TIP's for $\partial^+ \M$. Let $\gamma$ be any future-directed timelike
curve, and let $P= I^{-}(\gamma)$. If $\gamma$ is future inextendible, then it yields a TIP.  If $\gamma$ is continuously extendible to a point $p\in\M$,  then $P$ is equal to $I^-(p)$ and it is called a PIP, {\em proper indecomposable past set.} PIP's (and analogously PIF's) permit to identify $\M$ in the  future causal completion  $\M_c^+ = \M \cup \partial_c^+\M$ and, then, in the causal completion $\M_c = \M \cup \partial_c\M$.  
Parametrizing $\gamma$ with the $t$ coordinate of $\R\times M$, we have
$$\gamma(t)=(t,c(t)), \quad t\in [\alpha,
\Omega),  F(\dot c)  <1,$$
and, then,
$$\begin{array}{rl}
P  &= \{(t_{0},x_{0})\in \M : (t_{0},x_{0})\ll \gamma(t) \hbox{ for
	some }t\in [\alpha,\Omega) \} \\
&
= \{(t_{0},x_{0})\in \M : t_{0}<t-d^+(x_{0},c(t)) \hbox{ for some
} t\in [\alpha,\Omega)\} \\
&
= \{(t_{0},x_{0})\in \M : t_{0}<\lim_{t\rightarrow
	\Omega}(t-d^+(x_{0},c(t)))\}\\
&
=\{(t_{0},x_{0})\in V: t_{0}<b^+_{c}(x_{0})\}, 
\end{array}$$
where $b^+_{c}(x_{0}) = \lim_{t\rightarrow
	\Omega}(t-d^+(x_{0},c(t)))$  is the forward Busemann  function
of $c$ in $(M,d^+)$. So,  the set of  Busemann functions $B^+(M)$ for $(M,d^+)$ 
satisfies:\footnote{A subtlety is that the Busemann functions which can be constructed with the restriction $F(\dot c)\leq 1$ coincide with those constructed with $F(\dot c)< 1$.}

\begin{center}
	$\M_c^+(=\{$TIP's and PIP's on $\M \}$) $ \equiv$ $B^+(M)$   $\cup
	\{b_c\equiv \infty\}$, 
\end{center}
where the PIP's  correspond to converging $c$  (thus, necessarily, $\Omega
<\infty$)
and the TIP's to non-converging $c$, 
including the case $b_c\equiv \infty$ (this can be obtained with a curve  of  type $\gamma(t)=(t,x_0), t\geq 0$, for any $x_0\in M$). Notice that in the construction of $\M_c^+$
no quotient in the set of Busemann functions
is carried out.
Next, let us describe briefly $\partial_c \M$. We will restrict to its pointset  and chronological structures  and refer to  \cite[Chapter 5]{FHS_Memo} for the topological structure and full details.  

\subsubsection{The static case} 
\begin{figure}
\centering
	 \begin{tikzpicture}[scale=0.3]
	 \begin{scope}[shift={(3.4,0)}]
	 \draw[shift={(0 cm,5 cm)}]  (-3.2,1.3) -- (0.8,7);
	  \draw[shift={(0 cm,5 cm)}]  (4.7,2.1) -- (0.8,7);
	  \draw[purple] (0.8,-11) .. controls (1,-1) and (4,1) .. (0.8,12); 
	   \node at  (-2.1,8.9)   {$l^+$};
	      \node at  (-2,-8.9)   {$l^-$};
	   \node at  (1.5,12.9)   {$i^+$};
	     \node at  (1.5,-11.7)   {$i^-$};
	     \node at  (7,9.5)   {Future cone};
	      \node at  (7,6.5)   {$\partial_{\mathcal B}M$};
	        \node at  (7,-4.5)   {$\partial_{\mathcal B}M$};
	     \node at  (6.5,-8.9)   {Past cone};
	   \draw[purple] (0.8,-11) .. controls (-1,1) and (-4,1) .. (0.8,12); 
	\draw[dashed]  plot[smooth, tension=.7] coordinates {(-3.5,0.5) (-3,1.5) (-1,2.5) (1.5,2) (4,2.5) (5,1.5) (5,0.5) (2.5,-1) (0,-0.5) (-3,-1) (-3.5,0.5)};
		\draw[shift={(0 cm,5 cm)}]  plot[smooth, tension=.7] coordinates  {(-3.5,0.5) (-3,1.5) (-1,2.5) (1.5,2) (4,2.5) (5,1.5) (5,0.5) (2.5,-1) (0,-0.5) (-3,-1) (-3.5,0.5)};
		\draw[shift={(0 cm,-5 cm)}]  plot[smooth, tension=.7] coordinates  {(-3.5,0.5) (-3,1.5) (-1,2.5) (1.5,2) (4,2.5) (5,1.5) (5,0.5) (2.5,-1) (0,-0.5) (-3,-1) (-3.5,0.5)};
		\draw[shift={(0 cm,-5 cm)}]  (-3.3,-0.8) -- (0.8,-6);
	  \draw[shift={(0 cm,-5 cm)}]  (5.1,0.6) -- (0.8,-6);
	  \end{scope}
	\end{tikzpicture}

\caption{The causal boundary of a static spacetime. The  red lines joining $i_-$ and $i_+$ correspond to  points in $\partial_CM$.  
}
\label{FigStatic}
\end{figure}
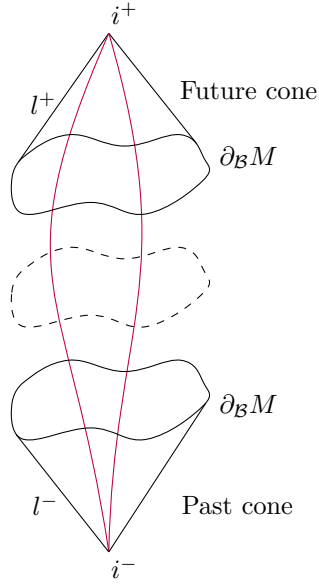

Recall that in the static case,  $\M $ can be regarded as a product $\R\times M$ with $F$ equal to a Riemannian metric by using conformal invariance.
Essentially, $\partial_c\M$ becomes a double cone with some lines connecting its apexes,  constructed as follows (see Figure \ref{FigStatic}):
\begin{enumerate}
	\item Two points $i^+=(M,\emptyset),i^-=(\emptyset,M)$ (which correspond to $b_c \equiv \infty$). They are the apexes of a symmetric double cone (invariant by $t\mapsto -t$) on $\partial_B M$.

	\item Two horismotic lines  (i.e., locally horismotic with no cut points) for each point of
	$\partial_\mathcal{B}M$, one of them $l^+$ ending at $i^+$
	and the other $l^-$ starting at $i^-$. This means that the line $l^+$ is  composed of boundary points of type  $(P,\emptyset)$ such that: (a) $l^+$ is totally ordered by the relation of inclusion for the first factor (i.e., $l^+$ is {\em locally horismotic}) ending at $P=\M$, and (b) no two points in $l^+$ are related by  the extended chronological relation $\overline{\ll}$ defined in Remark \ref{r_bordecausal} (2)  (i.e., $l^+$ has {\em no cut points}).
	
	\item A timelike line connecting $i^-,i^+$ for each point of the Cauchy boundary
	$\partial_C M$. Such a line is a continuous curve in $\M_c$ totally ordered by the extended chronological relation $\overline{\ll}$. 
	These are the only points in $\partial_c\M$ with non-trivial  pairings $(P,F)$.

\end{enumerate} 
Consistently with Remark \ref{r_bordecausal} (1), the spacetime is globally hyperbolic if the timelike lines do not exist, that is, if $\partial_CM=\emptyset$. 
%

\subsubsection{The general stationary case} 
For the sake of simplicity, we will assume that the causal completion $\M_c$ is {\em simple as a pointset}, that is, each TIP $P$ and each TIF $F$ determines unambiguously a  point of $\partial_c \M$ (thus, the situation in Figure \ref{Fig0b} cannot occur);  simple sufficient hypotheses to ensure this property are, for example, $M_B^\pm$ being Hausdorff or $d^+$ being extendible to  the  Cauchy boundary $\partial^+_CM$ as a generalized distance; see \cite[Figure 6.2]{FHS_Memo}.

Then, the three elements in the static picture must be modified allowing more generality, and  a fourth ingredient appears, namely, now:

\ben
\item $i^+=(M,\emptyset),i^-=(\emptyset,M)$ are regarded as apexes for two (non-symmetric) cones on $\partial_B^+M, \partial_B^-M$, respectively.

\item Each horismotic line $l^+$ (resp. $l^-$) appears for a point in 
$\partial_\mathcal{B}^+ M$ (resp. $\partial_\mathcal{B}^- M$).

\item  The timelike lines connecting $i^+$ and $i^-$ appear for each point in the Cauchy boundary for the symmetrized distance, $\partial_C^sM$.

\item The points in $\partial_C^\pm M\backslash \partial_C^sM$ determine  locally horismotic lines 
(eventually starting at $i^-$, ending at $i^+$, or both).

\een
Recall that the subtle last possibility cannot occur when  $\partial_C^\pm M = \partial_C^sM$, in particular, when $d^+$ extends to $\partial^+_CM$ as a generalized distance \cite[Proposition 3.28]{FHS_Memo}. 



\section{Lorentz-Finsler metrics and practical applications}\label{s6}
In this section we will take a look at some real-world situations where the interplay between Lorentz and Finsler geometries appears naturally. Specifically, we will focus on the wave propagation: Lorentz metrics and the spacetime viewpoint are essential when considering rheonomic (i.e., time-dependent) waves, whereas Finsler metrics effectively model the anisotropic (i.e., direction-dependent) case. The combination of both cases leads naturally to Lorentz-Finsler metrics. Interestingly, the applications can be generalized to any physical phenomenon that satisfies Huygens' principle, such as wildfires or seismic waves.

\subsection{Anisotropic wave propagation and Huygens' principle}
Recall from \S \ref{sec:time-dep} that the (time and direction-dependent) velocities of a moving object can be effectively described by the indicatrix $ \Sigma_t $ of a time-dependent Finsler metric $ F_t $ on $ \{t\}\times M $. Then the trajectories of the object are given in the (globally hyperbolic\footnote{This is not a restriction in any realistic situation (see \cite[Remark 3.2]{JPS}).}) spacetime $ \mathcal{M}=\mathds{R}\times M $ by the lightlike curves of the Lorentz-Finsler metric $ G=dt^2-F_t^2 $ introduced in \eqref{metricG},\footnote{Note that all the usual concepts about causality can be directly translated from the Lorentzian case to the Lorentz-Finsler metrics and the more general setting of cone structures (see for example \cite{JS20,Min}). However, due to the non-reversibility of the Finsler metrics, usually only future directions are considered.} since (considering $ t $-parametrized curves)
\begin{equation*}
\gamma(t)=(t,x(t)) \text{ lightlike} \Leftrightarrow G(\dot{\gamma}(t))=0 \Leftrightarrow x(t) \ F_t\text{-unit} \Leftrightarrow \dot{x}(t) \in \Sigma_t.
\end{equation*}

Consider now the following variation of Zermelo problem: instead of finding the fastest trajectory for a moving object between two prescribed points, we seek to determine the evolution of an anisotropic wave starting from $ \mathcal{S} \subset M $ with velocities given by $ \Sigma_t $. We will restrict ourselves to the ``mild wind'' case, i.e., $ \Sigma_t $ encloses the zero section,\footnote{Alternatively and more generally, we can consider $ \partial_t $ as an observer's vector field co-moving with the medium in which the wave propagates, as suggested at the end of \S \ref{s2}.  
} 
so that the wave propagates over $ M $ in all directions.

\begin{figure}
\centering
\includegraphics[width=0.9\textwidth]{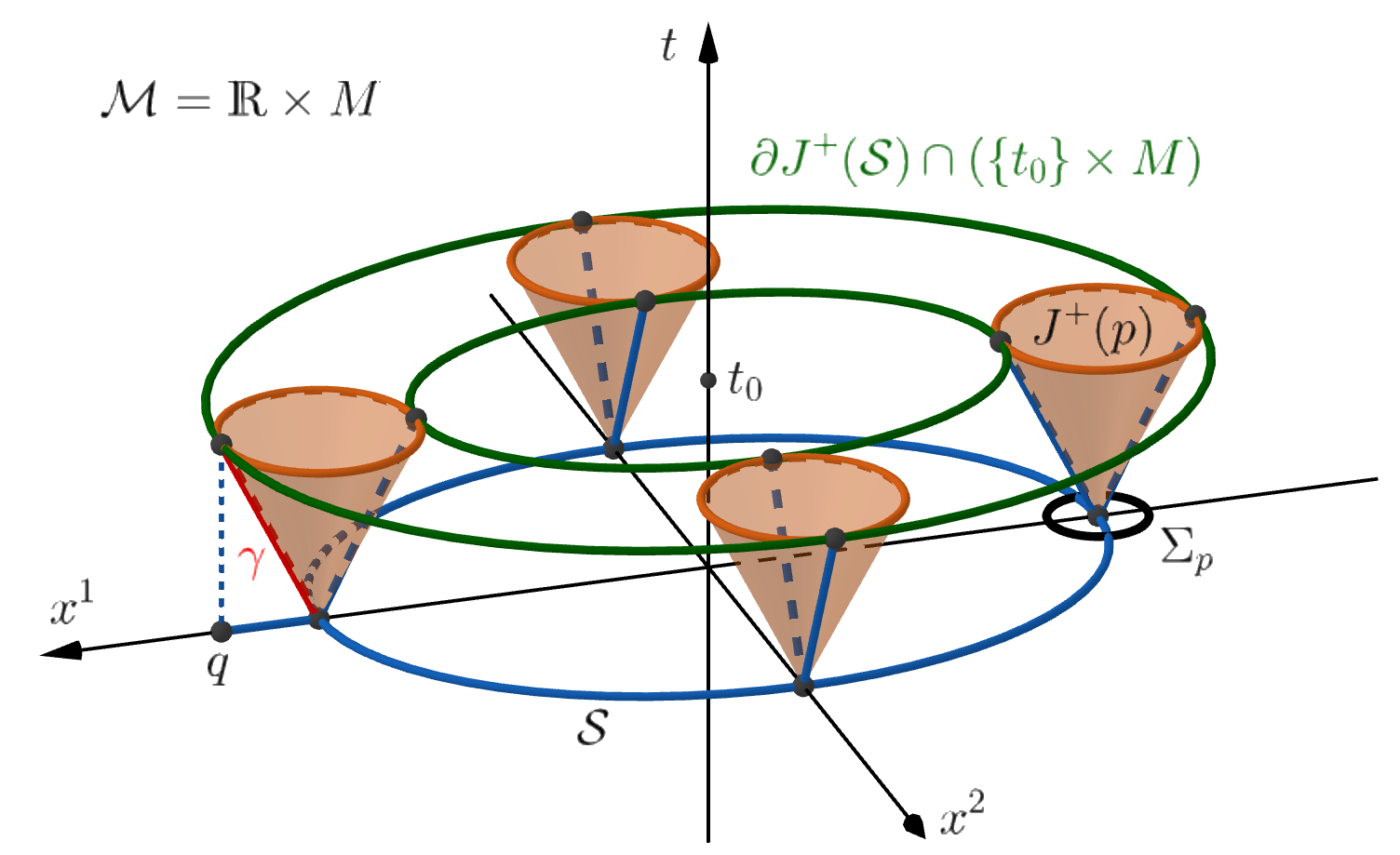}
\caption{A simple representation of the wave evolution in dimension $ n+1=3 $. The indicatrix $ \Sigma_p $ determines the velocity of the wave for each direction and the causal future $ J^+(p) $ shows the region in the spacetime that can be reached by the wave starting at $ p $. The envelope of all these chronological futures, $ \partial (\cup_{p \in \mathcal{S}}J^+(p)) = \partial J^+(\mathcal{S}) $, is generated by the first-arriving trajectories of the wave, so that the intersection with $ \{t_0\}\times M $ is the wavefront at any $ t_0 > 0 $ (when $ \mathcal{S} $ is a hypersurface of $ M $ there are always two wavefronts). A wave trajectory such as $ \gamma $ is a lightlike pregeodesic $ G $-orthogonal to $ \mathcal{S} $ and minimizes the propagation time from $ \mathcal{S} $ to $ q $.}
\label{fig:cones}
\end{figure}

Among all the possible trajectories of the wave (i.e., lightlike curves departing from $ \mathcal{S} $), we are interested in those that generate the wavefront. If $ \mathcal{S} = \{p\} $, namely, the wave starts from a single point $ p \in M $, then all the spacetime points that can be reached by the wave are given by the causal future $ J^+(p) $; accordingly, all the spatial points the wave passes through are the projection of $ J^+(p) $ on $ M $. Following this reasoning, $ \partial J^+(p) $ provides the outermost points reached by the wave and therefore, $ \partial J^+(p) \cap (\{t_0\}\times M) $ is the wavefront at each time $ t_0 > 0 $. Observe that $J^+(p)$ will be closed due to the global hyperbolicity of $ G $. 
 In the more general case where $ \mathcal{S} $ is a compact hypersurface of $ M $,\footnote{For simplicity, $ \mathcal{S} $ will be assumed to be a hypersurface of $ M $, although the results we present here can be generalized to any submanifold (see \cite{JPS}).} playing the role of the initial wavefront, we can apply Huygens' principle: each point of the front behaves as an independent source of the wave, and thus
\begin{equation*}
\text{Front}(t_1) = \partial\left( \cup_{p\in \text{Front}(t_0)} J^+(p) \right) \cap (\{t_1\}\times M) = \partial J^+(\text{Front}(t_0)) \cap (\{t_1\}\times M),
\end{equation*}
for any $t_1>t_0\geq 0$. In particular, if we put $ t_0=0 $, then the wavefront at any time $ t_1 > 0 $ is given by $ \partial J^+(\mathcal{S}) \cap (\{t_1\}\times M) $. So, the next step is to determine the trajectories that make up $ \partial J^+(\mathcal{S}) $ (see Figure \ref{fig:cones}).

\subsection{Solution in terms of Lorentz-Finsler geodesics}
So far, we know that a wave trajectory remaining in the wavefront must satisfy two conditions: it must be lightlike and it must be entirely contained in $ \partial J^+(\mathcal{S}) $. The next result states some crucial properties of this type of curves (see \cite[\S 4]{JPS} and \cite[\S 5]{JP}).

\begin{prop}\label{prop:time_min}
Let $ \gamma: [0,t_0] \rightarrow \mathcal{M} $ be a ($ t $-parametrized) causal curve entirely contained in $ \partial J^+(\mathcal{S}) $. Then:
\begin{enumerate}[(i)]
\item $ \gamma $ is a lightlike pregeodesic of $ (\mathcal{M},G) $ departing $ G $-orthogonally from $ \mathcal{S} $, i.e., $ g^G_{\dot{\gamma}(0)}(\dot{\gamma}(0),v)=0 $ for all $ v\in T_{\gamma(0)}\mathcal{S} $, being $ g^G $ the fundamental tensor of $ G $. In fact, $ \partial J^+(\mathcal{S}) $ admits a unique foliation by such pregeodesics.
\item $ \gamma $ is time-minimizing: for any $ p_0=(t_0,x_0) \in \textup{Im}(\gamma) $, $ \gamma $ is the first causal curve from $ \mathcal{S} $ to arrive at the vertical line $ t \mapsto (t,x_0) $.
\end{enumerate}
\end{prop}

Essentially, the wavefront is made up of lightlike pregeodesics that are time-minimizing in the sense that no other causal curve from $ \mathcal{S} $ arrives earlier at each of their points. In other words, these trajectories solve Zermelo problem from $ \mathcal{S} $ to any of their points.

Observe that at each point of $ \mathcal{S} $ there are exactly two lightlike $ G $-orthogonal directions: one points to the exterior of $ \mathcal{S} $ and the other to the interior. Therefore, there are two wavefronts: one heading outwards and the other inwards (see Figure \ref{fig:cones}). From now on we will focus on the one going outwards, since it is usually the most interesting from a practical viewpoint.\footnote{Nevertheless, the results we present here directly apply to the other wavefront simply by replacing ``outwards'' with ``inwards''.} Consider then the {\em wavemap}, defined as
\begin{equation*}
\begin{array}{rrll}
f \colon & [0,\infty)\times \mathcal{S} & \longrightarrow & \mathcal{M}\\
& (t,s) & \longmapsto & f(t,s)=(t,x(t,s)),
\end{array}
\end{equation*}
where, for each $ s \in \mathcal{S} $, $ t \mapsto f(t,s) $ is the $ t $-parametrized pregeodesic of $ G $ whose initial velocity is the unique lightlike vector $ G $-orthogonal to $ \mathcal{S} $ and pointing outwards. The curve $ t \mapsto f(t,s_0) $ represents the {\em spacetime trajectory} of the wave from $ s_0 \in \mathcal{S} $, being its projection $ t \mapsto x(t,s_0) $ the corresponding {\em spatial trajectory}.

Working in coordinates $ \{x^0 \coloneqq t,x^1,\ldots,x^n\} $, the following result characterizes the wavemap in terms of the ODE system for the $ t $-parametrized lightlike pregeodesics of $ G $ (see \cite[Theorem 4.11]{JPS}).

\begin{theorem}
For each $ s_0\in \mathcal{S} $, the wavemap $ f(t,s_0)=(t,x^1(t,s_0),\ldots,x^n(t,s_0)) $ is given by the following ODE system:
\begin{equation}
\label{eq:ode}
\ddot{x}^k = \sum_{i,j=0}^n \left(-\gamma_{\ ij}^k(\dot{f})\dot{x}^i\dot{x}^j + \gamma_{\ ij}^0(\dot{f})\dot{x}^i\dot{x}^j\dot{x}^k \right), \quad k=1,\ldots,n,
\end{equation}
along with the initial conditions:
\begin{itemize}
\item $ f(0,s_0) = s_0 \in \mathcal{S} $, and
\item $ \dot{f}(0,s_0) $ is lightlike, $ G $-orthogonal to $ \mathcal{S} $ and pointing outwards,
\end{itemize}
where $ \dot{f} = \dot{f}(t,s_0) = (1,\dot{x}^1(t,s_0),\ldots,\dot{x}^n(t,s_0)) $ denotes the velocity (tangent vector) of the curve $ t \mapsto f(t,s_0) $ and $ \gamma_{\ ij}^k $ are the formal Christoffel symbols of $ G $, defined as
\begin{equation*}
\gamma_{\ ij}^k(v) \coloneqq \frac{1}{2} \sum_{r=0}^n g^{kr}(v)\left(\frac{\partial g_{rj}}{\partial x^i}(v)+\frac{\partial g_{ri}}{\partial x^j}(v)-\frac{\partial g_{ij}}{\partial x^r}(v)\right), \quad i,j,k = 0,\ldots,n,
\end{equation*}
with $ g_{ij}(v) \coloneqq g^G_{v}(\partial_{x^i},\partial_{x^j}) $, for any lightlike $ v \in T\mathcal{M}$. 
\end{theorem}

\begin{remark}
In the time-independent anisotropic case (namely, $ F_t=F $ is time-independent), $ t $-parametrized lightlike pregeodesics of $ G $ project onto unit speed geodesics of $ F $ (recall Proposition \ref{lightlikegeo} for the Randers case). Namely, the spatial trajectories of the wave are $ F $-geodesics $ F $-orthogonal to $ \mathcal{S} $. Moreover, the propagation time is given by the $ F $-length: if $ \gamma: [0,t_0] \rightarrow M $ is a spatial trajectory of the wave, then $ t_0 $ coincides with the length of $ \gamma $ computed with $ F $. 
\end{remark}

\subsection{Cut points and determination of the wavefront}
As long as a spacetime trajectory remains in $ \partial J^+(\mathcal{S}) $, it keeps providing a point in the wavefront. However, this may not be the case for all $ t $. We define the {\em (null) cut function} as
\begin{equation*}
\begin{array}{rrll}
\textup{cut} \colon & \mathcal{S} & \longrightarrow & [0,\infty]\\
& s & \longmapsto & \textup{cut}(s)\coloneqq\text{Max}\{t:f(t,s)\in \partial J^+(\mathcal{S})\},
\end{array}
\end{equation*}
which also relates to the property of being time-minimizing. Specifically:
\begin{itemize}
\item If $ t \leq \textup{cut}(s_0) $, then $ f(t,s_0) \in \partial J^+(\mathcal{S}) $ is a point of the wavefront and the corresponding spacetime trajectory is time-minimizing.
\item If $ t > \textup{cut}(s_0) $, then $ f(t,s_0) \in I^+(\mathcal{S}) $ is not in the wavefront, so there exists another lightlike pregeodesic $ G $-orthogonal to $ \mathcal{S} $ and contained in $ \partial J^+(\mathcal{S}) $ that arrives earlier at the same spatial point $ x(t,s_0)\in M $ (recall that $ J^+(\mathcal{S}) $ is closed).
\end{itemize}
We call $ \textup{cut}(s) $ and $ (\textup{cut}(s),f(\textup{cut}(s),s)) $ the {\em cut instant} and {\em cut point} of the corresponding trajectory, respectively. Cut points are interesting from a practical viewpoint because they mark regions where different wave trajectories converge. More precisely (see \cite[Proposition A.1]{JPS2}):
\begin{prop}\label{prop:appendix}
Let $ (t_0,p_0)\in \mathcal{M} $ be the cut point of $ \gamma:t \mapsto f(t,s_0) $. Then, at least one of the following two properties holds: (a) $ (t_0,p_0) $ is the first intersection point of $ \gamma $ with another spacetime trajectory of the wave, or (b) $ (t_0,p_0) $ is the first focal point of $ \mathcal{S} $ along $ \gamma $ (it is possible for both conditions to hold simultaneously).
\end{prop}

In some situations, the detection of these points becomes crucial. In the case of wildfires, for example, cut points determine possible crossovers of the fire, which can become extremely dangerous both because of the increased heat intensity and because they may leave behind regions completely surrounded by the fire (see \cite[\S 4.2]{JPS2}).

Focusing now on the determination of the wavefront, observe that if $ t_0 < \textup{cut}(s) $ for all $ s \in \mathcal{S} $, then the curve $ s \mapsto f(t_0,s) $ exactly coincides with the wavefront at $ t =t_0 $, i.e., $ \text{Im}(f(t_0,s))=\partial J^+(\mathcal{S}) \cap (\{t_0\}\times M) $. This can be guaranteed at least for a small time, since there always exists some $\varepsilon>0$ such that $\textup{cut}(s)>\varepsilon$ for all $s\in \mathcal{S}$  (see \cite[Theorem 4.8]{JPS}). In general though, the wavefront will be given by all the trajectories which have not arrived yet at their cut points:
\begin{equation*}
\partial J^+(\mathcal{S}) = \{ f(t,s): t \leq \textup{cut}(s), s \in \mathcal{S} \}.
\end{equation*}
So, in order to obtain the wavefront we have to compute the wavemap through the ODE system \eqref{eq:ode}, but anytime a spacetime trajectory reaches its cut point it should be discarded, as subsequent points no longer lie on the wavefront. This is not demanding from a computational viewpoint, since each trajectory is independent from the others and those located beyond their cut points can be removed (or simply ignored) with no harm to the overall computation.

\subsection{The case of wildfires}
The theoretical setting we have presented in this section can be applied to model the propagation of any physical phenomena that behaves as a wave, i.e., that satisfies Huygens' principle. One of the most interesting examples is the case of wildfires (see \cite{M16,M17,JPS,JPS2}). Consider a wildfire that spreads over a surface $ \hat{M} \subset \mathds{R}^3 $. We can select a (global) coordinate chart $ (M,\hat{z}^{-1}) $, where $ \hat{z} $ is the graph
\begin{equation*}
\begin{array}{rrll}
\hat{z} \colon & M \subset \R^2 & \longrightarrow & \hat M \subset \R^3\\
& (x,y) & \longmapsto & \hat{z}(x,y) \coloneqq (x,y,z(x,y)),
\end{array}
\end{equation*}
and consider, as above, the spacetime $ (\mathcal{M},G=dt^2-F_t^2) $, where the indicatrix of $ F_t $ at each $ p=(t,x)\in \mathcal{M} $ provides the velocity of the fire for every direction. The Finslerian nature of the model in this specific case is absolutely essential, as several physical effects (mainly the slope and the wind) cause the propagation of the fire to be anisotropic.

If we know $ F_t $, i.e., we know the velocity of the fire at each point, time and direction, then solving \eqref{eq:ode} we obtain the evolution of the fire over time. The aim of a wildifire model is therefore to provide such $ F_t $.

In the isotropic case (without slope and wind), the indicatrix $ \Sigma_t $ is a sphere whose radius depends on the fuel and metereological conditions and may vary from one point to another (and over time) due to the change of vegetation, soil, moisture, etc. In order to include the isotropy caused by the wind, the most straightforward approximation is to consider that the wind displaces and deforms the sphere to an ellipse with a certain eccentricity depending on the wind strength. This approximation has been widely used since the experimental results by Anderson \cite{An} and the subsequent PDE system developed by Richards \cite{Rich} for the wavefront of a wildfire with an elliptical growth.\footnote{Richards' equations are equivalent to the ODE system \eqref{eq:ode} when $ \Sigma_t $ is an ellipse (see \cite[\S 5.2]{JPS}).} Richards' equations are still used nowadays by fire growth simulators such as FARSITE \cite{F} and Prometheus \cite{Prom}, which even extend the elliptical approximation to the isotropy caused by the slope, i.e., the wildfire becomes a displaced ellipse in the upward direction (since the fire moves faster upwards than downwards).

From a Finslerian viewpoint, the elliptical model translates into the metric $ F_t $ being of Randers type. Markvorsen was the first to propose the use of these metrics for wildfire modeling, transforming Richards' PDE into an ODE (the geodesic equations of the Randers metric; see \cite{M16}) and even developing a rheonomic Lagrangian viewpoint to include the time dependence (see \cite{M17}). For further developments using Randers metrics see \cite{D22}. Our work including Lorentz-Finsler metrics completes this theoretical framework and provides a full geometrical picture of the evolution of the wildfire in the most general situation. Specifically, there are mainly two important advantages of working with the Lorentz-Finsler setting over the classical elliptical one:
\begin{itemize}
\item {\em Flexibility}: the infinitesimal growth of the wildfire is not restricted to be elliptical and can adopt any other (strongly convex) pattern. In particular, the effect of the wind and the slope can be qualitatively different. For example, in \cite{JPS2} we have developed a specific model where the wind induces a sort of double semi-elliptical growth, which had already been pointed out by Anderson as the best experimental fitting in \cite{An}, while the slope generates the indicatrix of a reverse Matsumoto metric.\footnote{In its usual version, Matsumoto metrics effectively measure the  travel time for a walker on a slope, favoring the downward direction (see \cite{SS}).}
\item {\em Efficiency}: computationally speaking, solving an ODE is in general more efficient than solving a PDE. In addition, cut points represent a problem from the PDE viewpoint, since the firefront computed at an instant of time depends on the previously obtained, and therefore they must be corrected every time there is a crossover. This process has to be implemented through algorithms that are usually expensive in time and computing power (see for example \cite{F} and references therein). In comparison, this problem is greatly simplified in the ODE case, where we only need to remove the trajectories that reach their cut points, without even affecting the computation of the firefront (see \cite{JPS2}).
\end{itemize}



\section*{Acknowledgments}
MAJ and EPR were partially supported by the projects PGC2018-097046-B-I00 and PID2021-124157NB-I00, funded by MCIN/AEI/10.13039/501100011033/ "ERDF A way of making Europe", and also by Ayudas a proyectos para el desarrollo de investigaci\'{o}n cient\'{i}fica y t\'{e}cnica por grupos competitivos (Comunidad Aut\'{o}noma de la Regi\'{o}n de Murcia), included in the Programa Regional de Fomento de la Investigaci\'{o}n Cient\'{i}fica y T\'{e}cnica (Plan de Actuaci\'{o}n 2022) of the Fundaci\'{o}n S\'{e}neca-Agencia de Ciencia y Tecnolog\'{i}a de la Regi\'{o}n de Murcia, REF. 21899/PI/22. EPR and MS were partially supported by the project PID2020-116126GB-I00 funded by MCIN/AEI/10.13039/501100011033 and P20-01391 (PAIDI 2020, Junta de Andaluc\'{i}a), as well as the framework IMAG-Mar\'{i}a de Maeztu grant CEX2020-001105-M/AEI/10.13039/501100011033. EPR was also supported by Ayudas para la Formaci\'{o}n de Profesorado Universitario (FPU) from the Spanish Government.

%
%

\end{document}